\documentclass[a4paper,12pt]{article}

\usepackage[left=2cm,right=1.4cm, top=1.8cm,bottom=2.2cm,bindingoffset=0cm]{geometry}

\usepackage{verbatim}
\usepackage{amsmath}
\usepackage{amsthm}
\usepackage{amssymb}
\usepackage{delarray}
\usepackage{mathrsfs}
\usepackage{graphicx}
\usepackage{cite}

\begin{document}

\begin{center}
{\large\bf UNIFORM FULL STABILITY OF RECOVERING CONVOLUTIONAL PERTURBATION OF THE STURM--LIOUVILLE OPERATOR FROM THE SPECTRUM}
\\[0.5cm]

{\large\bf Sergey Buterin\footnote{Department of Mathematics, Saratov State University,
{\it email: buterinsa@info.sgu.ru}}} \\[0.2cm]
\end{center}

{\bf Abstract.} The perturbation of the Sturm--Liouville operator on a finite interval with Dirichlet boundary conditions by a convolution
operator is considered. Local stability and global unique solvability of the inverse problem of recovering the convolution kernel from the
spectrum, provided that the potential is given a priori, is known. In the present work, we establish uniform full stability of this inverse
problem involving a uniform estimate of deviations of the convolution kernel via deviations of the spectrum and the potential within balls
of any fixed radii.



{\bf \large 1. Introduction}
\\

The most complete results in the inverse spectral theory are known for differential operators (see, e.g., monographs
\cite{Mar77,Lev84,FY01,Yur02,Yur07,Yur00}). For integro-differential and other classes of nonlocal operators, the classical methods that
allow to obtain global solution of inverse problems such that the Gel'fand--Levitan method and the method of spectral mappings do not work.
Various aspects of inverse problems for integro-differential operators were studied in \cite{Mal, Yur84, Er, Yur91, But07, Kur, But10,
KurSh10, W,Yur14, BCh15,BS,BB17,Yur17-1, BondBut-18, But18,But18-2,Bon18-1,ButVas18,Ign18,Ign18-2,Bon18-2, Zol18, Bon19-1, Bon19-2,
Bon19-3,But19,But20,But07MM} and other works. One of the first substantial studies in this direction was undertaken in \cite{Yur91}, where
the boundary value problem ${\cal L}:={\cal L}(q,M)$ of the following form was considered:
\begin{equation}\label{equ}
-y''+q(x)y +\int_0^x M(x-t)y(t)\,dt=\lambda y, \;\; 0<x<\pi, \quad y(0)=y(\pi)=0,
\end{equation}
where $\lambda$ is the spectral parameter. We assume that $q(x)\in L_2(0,\pi)$ and $M(x)\in L_{2,\pi}$ being complex-valued functions,
where $L_{2,\pi}:=\{f(x):(\pi-x)f(x)\in L_2(0,\pi)\}.$ In \cite{Yur91}, it was established that the spectrum $\{\lambda_n\}_{n\ge1}$ of the
problem ${\cal L}$ has the form
\begin{equation}\label{asy}
\lambda_n=n^2+\omega+\varkappa_n,\quad \omega=\frac1\pi\int_0^\pi q(x)\,dx, \quad \{\varkappa_n\}\in l_2, \quad n\ge1,
\end{equation}
while the assumptions on the function $M(x)$ looked differently. Namely, it was assumed that
\begin{equation}\label{CondYur}
M_0(x):=(\pi-x)M(x), \, M_1(x):=\int_0^x M(t)\,dt\in L(0,\pi), \quad Q(x):=M_0(x)-M_1(x)\in L_2(0,\pi).
\end{equation}
However, according to Fubini's theorem, $M_0(x)\in L(0,\pi)$ implies $M_1(x)\in L(0,\pi).$ Moreover, from Lemma~2.4 in \cite{But07} (for
$\eta=1)$ it follows that (\ref{CondYur}) is equivalent to $M_0(x)\in L_2(0,\pi),$ i.e. $M(x)\in L_{2,\pi}$ (see also Remark~1 in
Section~6). Thus, we deal still with the same class of the functions $M(x)$ as in \cite{Yur91}, where the following inverse problem was
studied.

\medskip
{\bf Inverse Problem 1.} Given the spectrum $\{\lambda_n\}_{n\ge1},$ find the function $M(x),$ provided that the potential $q(x)$ is known
a priori.

\medskip
Developing the classical Borg method \cite{Bor} (see also \cite{BK19,FY01}), the uniqueness theorem for this inverse problem has been
proved and the following theorem was established, which gives its {\it local} solvability and stability (alternatively, see \cite{FY01}).

\medskip
{\bf Theorem 1. }{\it For any boundary value problem ${\cal L}={\cal L}(q,M)$ with the spectrum $\{\lambda_{n}\}_{n\ge 1}$ there exists
$\delta>0$ (which depends on ${\cal L}$) such that if an arbitrary sequence $\{\tilde\lambda_{n}\}_{n\ge 1}$ obeys the condition
$\|\{\lambda_n-\tilde\lambda_n\}\|_{l_2} \le\delta,$ then there exists a unique problem ${\cal L}(q,\tilde M)$ such that the sequence
$\{\tilde\lambda_{n}\}_{n\ge 1}$ is the spectrum of ${\cal L}(q,\tilde M).$ Moreover, the following estimates hold:
\begin{equation}\label{StabYur}
\|M_k-\tilde M_k\|_{L(0,\pi)}\le C_{\cal L}\|\{\lambda_n-\tilde\lambda_n\}\|_{l_2}, \quad k=0,1, \quad \|Q-\tilde Q\|_{L_{2}(0,\pi)}\le
C_{\cal L}\|\{\lambda_n-\tilde\lambda_n\}\|_{l_2},
\end{equation}
where $C_{\cal L}$ depends only on ${\cal L}.$}

\medskip
Note that the group of estimates (\ref{StabYur}) is equivalent to the estimate $\|M_0-\tilde M_0\|\le
C\|\{\lambda_n-\tilde\lambda_n\}\|_{l_2},$ i.e. $\|M-\tilde M\|_{2,\pi}\le C\|\{\lambda_n-\tilde\lambda_n\}\|_{l_2},$ where $C$ depends
only on $C_{\cal L}$ and vice versa (see Remark~1). Here and below, we use the designations $\|\cdot\|:=\|\cdot\|_{L_2(0,\pi)}$ and
$\|f\|_{2,\pi}:=\|(\pi-\,\cdot\,)f(\,\cdot\,)\|.$

In \cite{But10}, {\it global} solvability of Inverse Problem~1 was established, i.e. the following theorem holds.

\medskip
{\bf Theorem 2. }{\it Let a complex-valued function $q(x)\in L_2(0,\pi)$ be given. Then for any sequence of complex numbers
$\{\lambda_n\}_{n\ge1}$ of the form (\ref{asy}) there exists a unique (up to values on a set of measure zero) function $M(x)\in L_{2,\pi}$
such that $\{\lambda_n\}_{n\ge1}$ is the spectrum of the corresponding problem ${\mathcal L}(q,M).$}

\medskip
The proof was based on reducing the inverse problem to the so-called main nonlinear integral equation (see equation (\ref{ma_eq}) below)
and proving its global solvability in $L_{2,\pi}.$ This result was also briefly announced in \cite{But07MM}. Earlier, in \cite{But07}
Theorem~2 was established in the particular case when $q(x)\equiv const.$ We note that the case of arbitrary $q(x)\in L_2(0,\pi)$ is
essentially more difficult. Afterwards, evolvement of the approach suggested in \cite{But07} and \cite{But10} allowed to obtain global
solution also for other classes of integro-differential operators \cite{BCh15, BS, BB17, BondBut-18, But18, But18-2, Bon18-1, Ign18,
Ign18-2, Bon18-2, Bon19-1, Bon19-2, Bon19-3,But19}. In order to simplify proving solvability of the main equation in each new case, in
\cite{ButMal18} a general approach has been developed for solving nonlinear equations of this type and for proving their uniform stability.

In the present paper, by evolving the method used in \cite{But10} we obtain the following refinement of the stability part in Theorem~1.

\medskip
{\bf Theorem 3. }{\it For any fixed $r>0,$ there exists $C_r>0$ such that the estimate
\begin{equation}\label{UFStab}
\|M-\tilde M\|_{2,\pi}\le C_r(\|\{\varkappa_n-\tilde\varkappa_n\}\|_{l_2} +\|q-\tilde q\|)
\end{equation}
is fulfilled as soon as $\|\{\varkappa_n\}\|_{l_2}\le r,$ $\|\{\tilde\varkappa_n\}\|_{l_2}\le r$ and $\|q\|\le r,$ $\|\tilde q\|\le r.$
Here $\varkappa_n$ are determined in (\ref{asy}), while $\tilde\varkappa_n$ are determined by the analogous representation
$$
\tilde\lambda_n=n^2+\tilde\omega+\tilde\varkappa_n, \quad n\ge1,\quad \tilde\omega=\frac1\pi\int_0^\pi \tilde q(x)\,dx,
$$
where $\{\tilde\lambda_n\}_{n\ge1}$ is the spectrum of the boundary value problem $\tilde{\cal L}:={\cal L}(\tilde q,\tilde M).$}

\medskip
It is easy to see that estimate (\ref{UFStab}), in particular, implies estimates (\ref{StabYur}). However, unlike (\ref{StabYur}), estimate
(\ref{UFStab}) is uniform with respect to the spectra and the potentials of the both involved problems ${\cal L}$ and $\tilde{\cal L}.$
Moreover, since these problems are allowed to have different potentials, Theorem~2 actually gives uniform {\it full} stability of Inverse
Problem~1, i.e. uniform stability with respect to the complete set of input data. Note that the metric used in the right-hand side of
(\ref{UFStab}), unlike that in (\ref{StabYur}), admits different mean values $\omega$ and $\tilde\omega.$ Another advantage of the used
approach is insensitivity of the proof of Theorem~3 to the multiplicity of the spectrum, while in \cite{Yur91} for using a simpler analog
of Borg's method it was assumed that all eigenvalues were algebraically simple. Meanwhile, generalization of Borg's method to the case of
multiple spectrum appeared to be quite nontrivial and technical task (see \cite{BK19}). It should be mentioned that for the inverse problem
studied in \cite{But07} uniform stability was established in \cite{But20}. We note that uniform stability for the classical inverse
Sturm--Liouville problem in the selfadjoint case was obtained in \cite{SavShk} under the additional natural restriction preventing
neighboring eigenvalues to approach too close to each other. However, for Inverse Problem~1 as well as for the one in \cite{But07} no
similar restriction is necessary.

The paper is organized as follows. In the next section, we study the transformation operator related to the sine-type solution of equation
(\ref{equ}) as well as the characteristic function of the boundary value problem~${\cal L},$ and derive the main nonlinear integral
equation of the inverse problem. Section~3 is devoted to studying dependence of the transformation operator kernel on the functions $q(x)$
and $M(x).$ In Section~4, we provide the proof of global solvability and uniform full stability of the main equation. In Section~5, we
prove uniform stability of recovering the kernel of the characteristic function from its zeros. In Section~6, we give the proof of
Theorems~2 and~3, and provide a constructive procedure for solving the inverse problem (Algorithm~1).

Throughout the paper, one and the same symbol $C_r$ denotes {\it different} positive constants depending only on $r,$ and for any symbol
$\gamma$ we envisage the designation $\hat\gamma:=\gamma-\tilde\gamma.$
\\

{\large\bf  2. Characteristic function and transformation operator}
\\

Let $y=S(x,\lambda)=S(x,\lambda;q,M)$ be a solution of the equation in (\ref{equ}) under the initial conditions
\begin{equation}\label{ic}
S(0,\lambda)=0,\quad S'(0,\lambda)=1.
\end{equation}
Here and below, in order to emphasize dependence of a function $f(x_1,\ldots,x_n)$ on some functions $f_1,\ldots,f_m,$ sometimes we will
write $f(x_1,\ldots,x_n;f_1,\ldots,f_m).$

By virtue of uniqueness of the solution $S(x,\lambda),$ eigenvalues of ${\cal L}$ coincide with zeros (with account of multiplicity) of the
entire function
\begin{equation}\label{char}
\Delta(\lambda):=S(\pi,\lambda),
\end{equation}
which is called {\it characteristic function} of the problem ${\cal L}.$

For obtaining an appropriate representation of the function $S(x,\lambda),$ we need the following auxiliary assertion.

\medskip
{\bf Lemma 1. }{\it The integral equation
$$
F(x,t,\tau)=F_0(x,t,\tau)+\frac12\Big(\int_t^x q(s)\,ds\int_\tau^t F(s,\xi,\tau)\,d\xi+\int_\frac{t+\tau}{2}^t q(s)\,ds\int_\tau^{2s-t}\!
F(s,\xi,\tau)\,d\xi \qquad\qquad
$$
$$
\qquad\quad-\int_{\frac{\tau-t}{2}+x}^x q(s)\,ds\int_\tau^{2(s-x)+t} F(s,\xi,\tau)\,d\xi+\int_0^{t-\tau} M(s)\,ds\int_t^x
d\xi\int_\tau^{t-s} F(\xi-s,\eta,\tau)\,d\eta
$$
$$
+ \int_0^{t-\tau} M(s)\,ds\int_\frac{t+\tau+s}{2}^t d\xi\int_\tau^{2\xi-t-s} F(\xi-s,\eta,\tau)\,d\eta
$$
\begin{equation}\label{int_eq}
\quad\qquad\qquad-\int_0^{t-\tau} M(s)\,ds\int_{\frac{s+\tau-t}{2}+x}^x d\xi\int_\tau^{2(\xi-x)+t-s} F(\xi-s,\eta,\tau)\,d\eta\Big), \quad
0\le \tau\le t\le x\le \pi,
\end{equation}
with a continuous free term $F_0(x,t,\tau)$ has a unique solution $F(x,t,\tau)=F(x,t,\tau;q,M)$ being, in turn, a continuous function too.
Moreover, the following estimate holds:
\begin{equation}\label{int_eq_est}
|F(x,t,\tau)|\le F_0\exp(Ct), \quad 0\le\tau\le t\le x\le\pi,
\end{equation}
where
\begin{equation}\label{FC}
F_0=\max_{0\le\tau\le t\le x\le\pi}|F_0(x,t,\tau)|, \quad C=\int_0^\pi |q(s)|\,ds\ +\frac34\int_0^\pi (\pi-s)|M(s)|\,ds.
\end{equation}
}

\medskip
{\it Proof.} The method of successive approximations gives
\begin{equation}\label{msa}
F(x,t,\tau)=\sum_{k=0}^\infty F_k(x,t,\tau),
\end{equation}
where
$$
F_{k+1}(x,t,\tau)=\frac12\Big(\int_t^x q(s)\,ds\int_\tau^t F_k(s,\xi,\tau)\,d\xi+\int_\frac{t+\tau}{2}^t q(s)\,ds\int_\tau^{2s-t}
F_k(s,\xi,\tau)\,d\xi \qquad\qquad\qquad\qquad
$$
$$
\qquad -\int_{\frac{\tau-t}{2}+x}^x q(s)\,ds\int_\tau^{2(s-x)+t} F_k(s,\xi,\tau)\,d\xi +\int_0^{t-\tau} M(s)\,ds\int_t^x
d\xi\int_\tau^{t-s} F_k(\xi-s,\eta,\tau)\,d\eta
$$
$$
+\int_0^{t-\tau} M(s)\,ds\int_\frac{t+\tau+s}{2}^t d\xi\int_\tau^{2\xi-t-s} F_k(\xi-s,\eta,\tau)\,d\eta \qquad\qquad
$$
\begin{equation}\label{sa}
\qquad\qquad\qquad-\int_0^{t-\tau} M(s)\,ds\int_{\frac{s+\tau-t}{2}+x}^x d\xi\int_\tau^{2(\xi-x)+t-s} F_k(\xi-s,\eta,\tau)\,d\eta\Big),
\quad k\ge0.
\end{equation}
Let us show that
\begin{equation}\label{sa_est}
|F_k(x,t,\tau)|\le F_0\frac{(Ct)^k}{k!}, \quad 0\le\tau\le t\le x\le\pi, \quad k\ge0,
\end{equation}
where $F_0$ and $C$ are determined in (\ref{FC}). Indeed, for $k=0$ estimate (\ref{sa_est}) is obvious. Assuming it for some $k=j\ge0,$ one
can show it for $k=j+1.$ Indeed, according to (\ref{sa}), we have
$$
|F_{j+1}(x,t,\tau)|\le\frac12\Big(\int_\frac{t+\tau}{2}^x |q(s)|\,ds\int_\tau^t |F_j(s,\xi,\tau)|\,d\xi+\int_{\frac{\tau-t}{2}+x}^x
|q(s)|\,ds\int_\tau^t |F_j(s,\xi,\tau)|\,d\xi
$$
$$
+\int_0^{t-\tau} \!|M(s)|\,ds\int_\frac{t+\tau+s}{2}^x d\xi\int_\tau^t |F_j(\xi-s,\eta,\tau)|\,d\eta +\int_0^{t-\tau}
\!|M(s)|\,ds\int_{\frac{s+\tau-t}{2}+x}^x d\xi\int_\tau^t |F_j(\xi-s,\eta,\tau)|\,d\eta\Big).
$$
Since in the last two integrals $2x-t-\tau-s\le2(\pi-s)$ and $t-\tau-s\le\pi-s,$ respectively, substituting estimate (\ref{sa_est}) for
$k=j$ into the right-hand side of this inequality, we get the estimate
$$
|F_{j+1}(x,t,\tau)|\le F_0\frac{C^j}{2j!}\Big(\int_\frac{t+\tau}{2}^x |q(s)|\,ds\int_\tau^t \xi^j\,d\xi+\int_{\frac{\tau-t}{2}+x}^x
|q(s)|\,ds\int_\tau^t \xi^j\,d\xi \qquad\qquad\qquad\qquad
$$
$$
\qquad\qquad\qquad\qquad\qquad\qquad +\frac32\int_0^{t-\tau} (\pi-s)|M(s)|\,ds\int_\tau^t \eta^j\,d\eta\Big) \le
F_0\frac{(Ct)^{j+1}}{(j+1)!},
$$
which coincide with (\ref{sa_est}) for $k=j+1.$ Thus, the series in (\ref{msa}) converges uniformly in the pyramid $0\le \tau\le t\le x\le
\pi$ and, hence, its sum is a solution of equation (\ref{int_eq}).

It remains to show that $F_0(x,t,\tau)\equiv0$ implies $F(x,t,\tau)\equiv0.$ Indeed, assuming the zero free term we determine
$F_k(x,t,\tau)$ by formulae (\ref{sa}), having put $F_0(x,t,\tau):=F(x,t,\tau).$ Then, obviously, $F_k(x,t,\tau)=F(x,t,\tau)$ for all
$k\ge0$ and, by virtue of (\ref{sa_est}), we arrive at $F(x,t,\tau)\equiv0.$ $\hfill\Box$

\medskip
We note that the variable $\tau$ in equation (\ref{int_eq}) is, actually, a parameter, i.e. the assertion of Lemma~1 remains true, if one
fixes $\tau\in[0,\pi).$

The following lemma gives the transformation operator for the function $S(x,\lambda).$

\medskip
{\bf Lemma 2. }{\it Put $\rho^2=\lambda.$ The following representation holds:
\begin{equation}\label{op}
S(x,\lambda)=\frac{\sin\rho x}{\rho}+\int_0^x P(x,t)\frac{\sin\rho(x-t)}{\rho}\,dt, \quad 0\le x\le\pi,
\end{equation}
where the function
\begin{equation}\label{kop}
P(x,t)=P(x,t;q,M)=F(x,t,0;q,M)
\end{equation}
is the solution of equation (\ref{int_eq}) for $\tau=0$ and with the free term
\begin{equation}\label{ft}
F_0(x,t,0)=\frac12\Big(\int_{\frac{t}{2}}^{x-\frac{t}{2}}q(s)\,ds+\int_0^t (x-t)M(s)\,ds\Big).
\end{equation}
The function $P(x,t)$ is continuous in the triangle $0\le t\le x\le\pi.$ Moreover, $P(x,\,\cdot\,)\in W_2^1[0,x]$ for all $x\in(0,\pi]$ and
$P(\,\cdot\,,t)\in W_2^1[t,\pi]$ for all $t\in[0,\pi),$ and also
\begin{equation}\label{Kx0}
P(x,0)=\frac12\int_0^x q(t)\,dt, \quad P(x,x)=0, \quad 0\le x\le\pi.
\end{equation}}

{\it Proof.} By substitution it is easy to check that the Cauchy problem for the function $y=S(x,\lambda)$ consisting of the equation in
(\ref{equ}) and the initial conditions (\ref{ic}) is equivalent to the integral equation
\begin{equation}\label{eq_phi}
S(x,\lambda)=\frac{\sin\rho x}{\rho}+\int_0^x \frac{\sin\rho(x-t)}{\rho}\Big(q(t)S(t,\lambda) +\int_0^tM(t-s)S(s,\lambda)\,ds\Big)\,dt.
\end{equation}
Substituting the claimed representation (\ref{op}) into equation (\ref{eq_phi}) and multiplying by $\rho,$ we obtain
\begin{equation}\label{K_rho}
\int_0^x P(x,t)\sin\rho(x-t)\,dt=\sum_{\nu=1}^4 {\mathcal P}_\nu(x,\lambda),
\end{equation}
where
$$
{\mathcal P}_1(x,\lambda)=\int_0^x \sin\rho(x-t)q(t)\,dt\int_0^t\cos\rho s\,ds,
$$
$$
{\mathcal P}_2(x,\lambda)=\int_0^x \sin\rho(x-t)\,dt \int_0^tM(t-s)\,ds \int_0^s\cos\rho \xi\,d\xi,
$$
$$
{\mathcal P}_3(x,\lambda)=\int_0^x \sin\rho(x-t)q(t)\,dt\int_0^t P(t,t-s)\,ds\int_0^s\cos\rho \xi\,d\xi,
$$
$$
{\mathcal P}_4(x,\lambda)=\int_0^x \sin\rho(x-t)\,dt \int_0^tM(t-s)\,ds\int_0^s P(s,s-\xi)\,d\xi\int_0^\xi\cos\rho\eta\,d\eta.
$$
Since
$$
\sin\rho(x-t)\cos\rho s=\frac12\Big(\sin\rho(x-t+s)+\sin\rho(x-t-s)\Big),
$$
changing the variables and the order of integration, we get
\begin{equation}\label{2eq2.9.3}
{\mathcal P}_1(x,\lambda)=\frac12\int_0^x q(t)\,dt\int_{x-2t}^x\sin\rho s\,ds =\frac12\int_0^x \sin\rho(x-t) \,dt
\int_\frac{t}2^{x-\frac{t}2}q(s)\,ds,
\end{equation}
$$
{\mathcal P}_2(x,\lambda)=\frac12\int_0^x \Big(\int_{x-t}^x\sin\rho s\,ds\int_0^{x-s}M(\xi)\,d\xi + \int_{x-2t}^{x-t}\sin\rho s\,ds
\int_0^{2t-x+s}M(\xi)\,d\xi\Big)dt \qquad\qquad\qquad\qquad
$$
\begin{equation}\label{2eq2.9.4} =\frac12\int_0^x\sin\rho t\,dt
\int_{x-t}^x\,ds\int_0^{x-t}M(\xi)\,d\xi=\frac12\int_0^x(x-t)\sin\rho(x-t)\,dt \int_0^t M(s)\,ds,
\end{equation}
$$
{\mathcal P}_3(x,\lambda)=\frac12\int_0^x q(t) \Big(\int_{x-t}^x\sin\rho s\,ds\int_{s-x+t}^tP(t,t-\xi)\,d\xi +\int_{x-2t}^{x-t}\sin\rho
s\,ds \int_{x-t-s}^tP(t,t-\xi)\,d\xi\Big)dt\qquad\qquad\qquad\quad
$$
$$
=\frac12\int_0^x\sin\rho(x-t)\Big(\int_t^xq(s)\,ds\int_0^tP(s,\xi)\,d\xi +
\int_\frac{t}2^tq(s)\,ds\int_0^{2s-t}P(s,\xi)\,d\xi\quad\quad\quad
$$
\begin{equation}\label{2eq2.9.5}
\;\quad\qquad\qquad\qquad\qquad\qquad\qquad\qquad\qquad\qquad\qquad\qquad-\int_{x-\frac{t}2}^xq(s)\,ds\int_0^{2(s-x)+t}P(s,\xi)\,d\xi\Big)dt,
\end{equation}
$$
{\mathcal P}_4(x,\lambda)=\frac12\int_0^x \Big(\int_{x-t}^x\sin\rho s\,ds\int_0^{x-s}M(\xi)\,d\xi\int_0^{x-s-\xi}P(t-\xi,\eta)\,d\eta
\qquad\qquad\qquad\qquad\qquad\quad\quad\quad
$$
$$
+\int_{x-2t}^{x-t}\sin\rho s\,ds \int_0^{2t-x+s}M(\xi)\,d\xi \int_0^{2t-x+s-\xi}P(t-\xi,\eta)\,d\eta\Big)dt \qquad\qquad
$$
$$
=\frac12\int_0^x\sin\rho(x-t)\Big(\int_0^tM(s)\,ds\int_t^x d\xi \int_0^{t-s}P(\xi-s,\eta)\,d\eta
$$
$$
\qquad +\int_0^tM(s)\,ds\int_{\frac{t+s}2}^t d\xi\int_0^{2\xi-t-s}P(\xi-s,\eta)\,d\eta
$$
\begin{equation}\label{2eq2.9.6}
\qquad\qquad\qquad\qquad\qquad\qquad\qquad\qquad-\int_0^tM(s)\,ds\int_{\frac{s-t}2+x}^x
d\xi\int_0^{2(\xi-x)+t-s}P(\xi-s,\eta)\,d\eta\Big)dt.
\end{equation}
According to (\ref{2eq2.9.3})--(\ref{2eq2.9.6}), identity (\ref{K_rho}) holds for all $\rho\in{\mathbb C}$ if and only if the function
$P(x,t)$ satisfies the linear integral equation
$$
P(x,t)=\frac12\Big(\int_\frac{t}2^{x-\frac{t}2}q(s)\,ds +\int_0^t (x-t)M(s)\,ds +\int_t^xq(s)\,ds\int_0^tP(s,\xi)\,d\xi
\qquad\qquad\qquad\qquad\qquad
$$
$$
+\int_\frac{t}2^tq(s)\,ds\int_0^{2s-t}P(s,\xi)\,d\xi -\int_{x-\frac{t}2}^xq(s)\,ds\int_0^{2(s-x)+t}P(s,\xi)\,d\xi \qquad\qquad\qquad\qquad
$$
$$
\quad+\int_0^tM(s)\,ds\int_t^x d\xi \int_0^{t-s}P(\xi-s,\eta)\,d\eta+ \int_0^tM(s)\,ds\int_{\frac{t+s}2}^t
d\xi\int_0^{2\xi-t-s}P(\xi-s,\eta)\,d\eta
$$
\begin{equation}\label{int_eq1}
\qquad\qquad\qquad\qquad\qquad -\int_0^tM(s)\,ds\int_{\frac{s-t}2+x}^x d\xi\int_0^{2(\xi-x)+t-s}P(\xi-s,\eta)\,d\eta\Big), \quad 0\le t\le
x\le \pi,
\end{equation}
which, in turn, in accordance with (\ref{kop}) and (\ref{ft}), is equivalent to  equation (\ref{int_eq}) for $\tau=0.$

The rest properties of the kernel $P(x,t)$ immediately follow from the form of equation (\ref{int_eq1}). $\hfill\Box$

\medskip
Using (\ref{char}) and Lemma~2, by simple calculations we arrive at the following lemma, which gives a fundamental representation of the
characteristic function.

\medskip
{\bf Lemma 3. }{\it The characteristic function of the problem ${\mathcal L}$ has the form
\begin{equation}\label{char1}
\Delta(\lambda)=\frac{\sin\rho\pi}{\rho}-\omega\pi\frac{\cos\rho\pi}{2\rho^2} +\int_0^\pi v(x)\frac{\cos\rho x}{\rho^2}\,dx, \quad v(x)\in
L_2(0,\pi).
\end{equation}
Here $\omega$ is determined in (\ref{asy}) and
\begin{equation}\label{ma_eq}
-v(\pi-x)=R(\pi,x;q,M), \quad 0<x<\pi,
\end{equation}
where
\begin{equation}\label{K_12}
R(x,t;q,M)=\frac{\partial}{\partial t}P(x,t;q,M).
\end{equation}}

Taking entireness of the function $\Delta(\lambda)$ into account, it is easy to see that
\begin{equation}\label{means}
\int_0^\pi v(x)\,dx=\frac{\omega\pi}2=\frac12\int_0^\pi q(x)\,dx.
\end{equation}

We remind that asymptotics (\ref{asy}) can be established using representation (\ref{char1}) by the known method (see, e.g., \cite{Mar77})
involving Rouch\'e's theorem. Moreover, using Hadamard's factorization theorem, by the standard approach (see, e.g., \cite{FY01}) one can
prove the following lemma.

\medskip
{\bf Lemma 4. }{\it Any function of the form (\ref{char1}) is uniquely determined by its zeros $\lambda_n,$ $n\ge1.$ Moreover, the
following representation holds:
\begin{equation}\label{char2}
\Delta(\lambda)=\pi\prod_{n=1}^\infty \frac{\lambda_n-\lambda}{n^2}.
\end{equation}}

\medskip
The following lemma gives the inverse assertion (see, e.g., Lemma~3.3 in \cite{But07}).

\medskip
{\bf Lemma 5. }{\it Let an arbitrary complex sequence $\{\lambda_n\}_{n\ge1}$ of the form (\ref{asy}) be given. Then the function
$\Delta(\lambda)$ determined by (\ref{char2}) has the form (\ref{char1}) with a certain function $v(x)\in L_2(0,\pi)$
obeying~(\ref{means}).}

\medskip
Relation (\ref{ma_eq}) can be considered as a nonlinear equation with respect to the function $M(x),$ to which we refer as {\it main
nonlinear equation} or shortly {\it main equation} of the inverse problem. Let an arbitrary function $q(x)\in L_2(0,\pi)$ be given. Up to
now we established that for any function $M(x)\in L_{2,\pi}$ the function $v(x)$ determined by formula (\ref{ma_eq}) belongs to
$L_2(0,\pi)$ and obeys the condition (\ref{means}). The following theorem gives the inverse assertion, which occupies the central place in
our approach to solving Inverse Problem~1.

\medskip
{\bf Theorem 4. }{\it (i) For any complex-valued functions $q(x),\,v(x)\in L_2(0,\pi),$ satisfying relation (\ref{means}), the nonlinear
equation (\ref{ma_eq}) has a unique solution $M(x)\in L_{2,\pi}.$

(ii) Fix $r>0$ and, besides (\ref{means}), let $\|v\|\le r$ and $\|q\|\le r.$ Then the estimate $\|M\|_{2,\pi}\le C_r$ is fulfilled. }

\medskip
The proof of the first part of Theorem~4, which gives global solvability of the main nonlinear equation (\ref{ma_eq}), was briefly given in
\cite{But10}. Here we  provide the complete detailed proof of this theorem including its second part, which plays an important role in
proving uniform full stability of the main equation (see Theorem~5 in Section~4).

The proof of Theorem~4 is based on further properties of the kernel $P(x,t;q,M)$ as an operator on the pair of functions $q(x)$ and $M(x),$
which will be established in the next section.
\\

{\large\bf 3. Properties of the mapping $P(x,t;\,\cdot\,,\,\cdot\,)$}
\\

In the present section, we reveal some important properties of the kernel $P(x,t;q,M)$ regarding its dependence on the functions $q(x)$ and
$M(x),$ which allow to prove global solvability of the main nonlinear equation (\ref{ma_eq}) and its uniform full stability. First, we make
the following observation.

\medskip
{\bf Observation 1.} For each fixed $\delta\in(0,\pi]$ the linear integral equation (\ref{int_eq}) can be restricted to the set
$$
{\mathcal D}_\delta:=\Big\{(x,t,\tau):0\le\tau\le t\le\min\{\delta,x\},\; x\le \pi\Big\}.
$$
In other words, for $(x,t,\tau)\in{\mathcal D}_\delta$ the right-hand side of equation (\ref{int_eq}) depends on values of the unknown
function $F(x,t,\tau)=F(x,t,\tau;q,M)$ only on the subset ${\mathcal D}_\delta.$ Obviously, the solution of the restricted equation
coincides with the restriction to ${\mathcal D}_\delta$ of the solution of the initial one. Hence, the function $F(x,t,\tau;q,M)$ on
${\mathcal D}_\delta$ depends on values of the function $M(s)$ only on the interval $(0,\delta).$ In particular, due to (\ref{kop}) and
(\ref{ft}) or (\ref{int_eq1}), the kernel $P(x,t)=P(x,t;q,M)$ on the trapezium
$$
D_\delta:=\Big\{(x,t):0\le t\le\min\{\delta,x\},\; x\le \pi\Big\}
$$
depends on values of $M(s)$ only for $s\in(0,\delta).$

\medskip
{\bf Lemma 6. }{\it The following representation holds:
\begin{equation}\label{3.4-1}
P_M(x,t):= P(x,t;q,M)-P(x,t;q,\tilde M) =\int_0^t F(x,t,\tau;q,M,\tilde M)\hat M(\tau)\,d\tau, \quad 0\le t\le x\le\pi,
\end{equation}
where the function $F(x,t,\tau)=F(x,t,\tau;q,M,\tilde M)$ is a solution of equation (\ref{int_eq}) with the free term
$$
F_0(x,t,\tau)=F_0(x,t,\tau;q,\tilde M)=\frac12\Big(x-t+\int_t^x ds\int_0^{t-\tau} P(s-\tau,\xi;q,\tilde
M)\,d\xi\qquad\qquad\qquad\qquad\qquad
$$
\begin{equation}\label{3.5-1}
\;\;\quad\quad\qquad+\int_\frac{t+\tau}{2}^t ds\int_0^{2s-t-\tau} P(s-\tau,\xi;q,\tilde M)\,d\xi -\int_{\frac{\tau-t}{2}+x}^x
ds\int_0^{2(s-x)+t-\tau} P(s-\tau,\xi;q,\tilde M)\,d\xi\Big).
\end{equation}}

\medskip
{\it Proof.} Termwise subtracting equation (\ref{int_eq1}) for the function $P(x,t;q,\tilde M),$ i.e. after substituting $\tilde M(x)$
instead of $M(x),$ from the unchanged equation (\ref{int_eq1}), we arrive at the following integral equation with respect to the function
$P_M(x,t):$
$$
P_M(x,t)=\int_0^t F_0(x,t,\tau;q,\tilde M) \hat M(\tau)\,d\tau +\frac12\Big(\int_t^xq(s)\,ds\int_0^tP_M(s,\xi)\,d\xi
\qquad\qquad\qquad\qquad\qquad\;\;
$$
$$
+\int_\frac{t}2^tq(s)\,ds\int_0^{2s-t}P_M(s,\xi)\,d\xi -\int_{x-\frac{t}2}^xq(s)\,ds\int_0^{2(s-x)+t}P_M(s,\xi)\,d\xi \qquad\qquad\qquad
$$
$$
\quad\qquad+\int_0^tM(s)\,ds\int_t^x d\xi \int_0^{t-s}P_M(\xi-s,\eta)\,d\eta+ \int_0^tM(s)\,ds\int_{\frac{t+s}2}^t
d\xi\int_0^{2\xi-t-s}P_M(\xi-s,\eta)\,d\eta
$$
\begin{equation}\label{int_eq2}
\quad\qquad\qquad\qquad\qquad -\int_0^tM(s)\,ds\int_{\frac{s-t}2+x}^x d\xi\int_0^{2(\xi-x)+t-s}P_M(\xi-s,\eta)\,d\eta\Big), \quad 0\le t\le
x\le \pi,
\end{equation}
where the function $F_0(x,t,\tau;q,\tilde M)$ is determined by (\ref{3.5-1}). Substituting the target representation (\ref{3.4-1}) into
equation (\ref{int_eq2}), we conclude that the right-hand side of (\ref{3.4-1}) is a solution of (\ref{int_eq2}) if and only if the
following relation holds:
\begin{equation}\label{3.6}
\int_0^t F(x,t,\tau;q,M,\tilde M)\hat M(\tau)\,d\tau=\int_0^t F_0(x,t,\tau;q,\tilde M)\hat M(\tau)\,d\tau  + \frac12\sum_{k=1}^3({\mathcal
Q}_k(x,t)+{\mathcal M}_k(x,t)),
\end{equation}
where
$$
{\mathcal Q}_1(x,t):= \int_t^x q(s)\,ds\int_0^t\,d\xi\int_0^\xi F(s,\xi,\tau;q,M,\tilde M)\hat M(\tau)\,d\tau
\qquad\qquad\qquad\qquad\qquad\qquad\qquad\quad
$$
$$
\;= \int_t^x q(s)\,ds\int_0^t \hat M(\tau)\,d\tau \int_\tau^t F(s,\xi,\tau;q,M,\tilde M)\,d\xi
$$
\begin{equation}\label{calQ_1}
\qquad\qquad \qquad\qquad\qquad\qquad\qquad =\int_0^t \hat M(\tau) \,d\tau \int_t^x q(s)\,ds\int_\tau^t F(s,\xi,\tau;q,M,\tilde M) \,d\xi,
\end{equation}
$$
{\mathcal Q}_2(x,t):= \displaystyle \int_\frac{t}{2}^t q(s)\,ds\int_0^{2s-t}\,d\xi\int_0^\xi
F(s,\xi,\tau;q,M,\tilde M)\hat M(\tau)\,d\tau \qquad\qquad\qquad\qquad\qquad\qquad\quad\quad
$$
$$
=\int_\frac{t}{2}^t q(s)\,ds\int_0^{2s-t}\hat M(\tau)\,d\tau\int_\tau^{2s-t} F(s,\xi,\tau;q,M,\tilde M)\,d\xi
$$
\begin{equation}\label{calQ_2}
\qquad\qquad\qquad \qquad\qquad\qquad = \int_0^t \hat M(\tau)\,d\tau \int_\frac{t+\tau}{2}^t q(s)\,ds\int_\tau^{2s-t}
F(s,\xi,\tau;q,M,\tilde M)\,d\xi,
\end{equation}
$$
{\mathcal Q}_3(x,t):= \displaystyle -\int_{x-\frac{t}{2}}^x q(s)\,ds\int_0^{2(s-x)+t}\,d\xi\int_0^\xi F(s,\xi,\tau;q,M,\tilde M)\hat
M(\tau)\,d\tau\qquad\qquad\qquad\qquad\quad\quad\;
$$
$$
\qquad= -\int_{x-\frac{t}{2}}^x q(s)\,ds\int_0^{2(s-x)+t}\hat M(\tau)\,d\tau\int_\tau^{2(s-x)+t}
F(s,\xi,\tau;q,M,\tilde M)\,d\xi
$$
\begin{equation}\label{calQ_3}
\qquad\qquad \qquad\quad\quad =-\int_0^t \hat M(\tau)\,d\tau \int_{\frac{\tau-t}{2}+x}^x q(s)\,ds\int_\tau^{2(s-x)+t}
F(s,\xi,\tau;q,M,\tilde M)\,d\xi,
\end{equation}
$$
{\mathcal M}_1(x,t):=\int_0^t M(s)\,ds\int_t^x d\xi\int_0^{t-s}d\eta\int_0^\eta F(\xi-s,\eta,\tau;q,M,\tilde M)\hat M(\tau)\,d\tau
\qquad\qquad\qquad\qquad
$$
$$
=\int_0^t M(s)\,ds\int_0^{t-s}\hat M(\tau)\,d\tau\int_t^x d\xi\int_\tau^{t-s} F(\xi-s,\eta,\tau;q,M,\tilde M)\,d\eta
$$
\begin{equation}\label{calM_1}
\quad\qquad\qquad\qquad\qquad=\int_0^t \hat M(\tau)\,d\tau\int_0^{t-\tau}M(s)\,ds\int_t^x d\xi\int_\tau^{t-s} F(\xi-s,\eta,\tau;q,M,\tilde
M)\,d\eta,
\end{equation}
$$
{\mathcal M}_2(x,t):=\int_0^t M(s)\,ds\int_\frac{t+s}{2}^t d\xi\int_0^{2\xi-t-s}\,d\eta\int_0^\eta F(\xi-s,\eta,\tau;q,M,\tilde M)\hat
M(\tau)\,d\tau \qquad\qquad\qquad\;
$$
$$
\qquad\qquad\qquad=\int_0^t M(s)\,ds\int_\frac{t+s}{2}^t d\xi\int_0^{2\xi-t-s}\hat M(\tau)\,d\tau\int_\tau^{2\xi-t-s}
F(\xi-s,\eta,\tau;q,M,\tilde M)\,d\eta \qquad\qquad\quad\;
$$
$$
\qquad=\int_0^t M(s)\,ds\int_0^{t-s}\hat M(\tau)\,d\tau\int_\frac{t+\tau+s}{2}^t d\xi\int_\tau^{2\xi-t-s} F(\xi-s,\eta,\tau;q,M,\tilde
M)\,d\eta
$$
\begin{equation}\label{calM_2}
\qquad\qquad\qquad=\int_0^t \hat M(\tau)\,d\tau\int_0^{t-\tau}M(s)\,ds\int_\frac{t+\tau+s}{2}^t d\xi\int_\tau^{2\xi-t-s}
F(\xi-s,\eta,\tau;q,M,\tilde M)\,d\eta,
\end{equation}
$$
{\mathcal M}_3(x,t):=-\int_0^t M(s)\,ds\int_{\frac{s-t}{2}+x}^x d\xi\int_0^{2(\xi-x)+t-s}\,d\eta\int_0^\eta F(\xi-s,\eta,\tau;q,M,\tilde
M)\hat M(\tau)\,d\tau\qquad\quad
$$
$$
\quad\quad=-\int_0^t M(s)\,ds\int_{\frac{s-t}{2}+x}^x d\xi\int_0^{2(\xi-x)+t-s}\hat M(\tau)\,d\tau\int_\tau^{2(\xi-x)+t-s}
F(\xi-s,\eta,\tau;q,M,\tilde M)\,d\eta\qquad
$$
$$
=-\int_0^t M(s)\,ds\int_0^{t-s}\hat M(\tau)\,d\tau \int_{\frac{s+\tau-t}{2}+x}^x d\xi\int_\tau^{2(\xi-x)+t-s} F(\xi-s,\eta,\tau;q,M,\tilde
M)\,d\eta
$$
\begin{equation}\label{calM_3}
\quad\quad=-\int_0^t \hat M(\tau)\,d\tau\int_0^{t-\tau}M(s)\,ds\int_{\frac{s+\tau-t}{2}+x}^x d\xi\int_\tau^{2(\xi-x)+t-s}
F(\xi-s,\eta,\tau;q,M,\tilde M)\,d\eta.
\end{equation}
Taking (\ref{calQ_1})--(\ref{calM_3}) into account, we conclude that if the function $F(x,t,\tau;q,M,\tilde M)$ obeys conditions of the
lemma, then equality (\ref{3.6}) is fulfilled. Thus, the both sides of (\ref{3.4-1}) satisfy one and the same equation (\ref{int_eq2}),
having a unique solution, which finishes the proof. $\hfill\Box$

\medskip
This lemma will be used for proving stability of the main equation, while the following its corollary, giving stepwise linearizability of
the operator $P(x,t;q,\,\cdot\,),$ plays a crucial role in proving Theorem~4. In what follows, for any fixed $\delta\in(0,\pi/2]$ we use
the designations
$$
M_1(x):=\left\{\begin{array}{cl}M(x), & x\in(0,\delta),\\[2mm] 0, & x\in(\delta,2\delta), \end{array}\right. \quad
M_2(x):=\left\{\begin{array}{cl}0, & x\in(0,\delta),\\[2mm] M(x), & x\in(\delta,2\delta). \end{array}\right.
$$

\medskip
{\bf Corollary 1. }{\it For any $\delta\in(0,\pi/2]$ the representation
\begin{equation}\label{2eq3.4}
P(x,t;q,M)=P(x,t;q,M_1)+\int_0^t F(x,t,\tau;q,M_1,M_1)M_2(\tau)\,d\tau, \quad (x,t)\in D_{2\delta},
\end{equation}
holds, where the function $F(x,t,\tau;q,M,\tilde M)$ is determined in Lemma~6.}

\medskip
{\it Proof.} According to Observation~1, formula (\ref{3.4-1}) for $(x,t)\in D_{2\delta}$ involves $M(x)$ and $\tilde M(x)$ only for $x\in
(0,2\delta).$ Taking $\tilde M(x)=M_1(x),$ we get $\hat M(x)=M_2(x)$ and, hence, (\ref{3.4-1}) takes the form
$$
P(x,t;q,M)=P(x,t;q,M_1)+\int_0^t F(x,t,\tau;q,M,M_1)M_2(\tau)\,d\tau, \quad (x,t)\in D_{2\delta},
$$
where it remains to show that $F(x,t,\tau;q,M,M_1)=F(x,t,\tau;q,M_1,M_1)$ for $(x,t,\tau)\in {\mathcal D}_{2\delta}.$ Indeed, since $\hat
M(x)=M_2(x)$ and
$$
\int_0^t M_2(\tau)\,d\tau\int_0^{t-\tau}M_2(s)\,ds\int_t^x d\xi\int_\tau^{t-s} F(\xi-s,\eta,\tau;q,M,\tilde M)\,d\eta=0, \quad 0\le
t\le2\delta,
$$
formula (\ref{calM_1}) takes the form
\begin{equation}\label{calM_1-1}
{\mathcal M}_1(x,t)=\int_0^t M_2(\tau)\,d\tau\int_0^{t-\tau}M_1(s)\,ds\int_t^x d\xi\int_\tau^{t-s} F(\xi-s,\eta,\tau;q,M,\tilde M)\,d\eta.
\end{equation}
Analogously one can get
\begin{equation}\label{calM_2-1}
{\mathcal M}_2(x,t)=\int_0^t M_2(\tau)\,d\tau\int_0^{t-\tau}M_1(s)\,ds\int_\frac{t+\tau+s}{2}^t d\xi\int_\tau^{2\xi-t-s}
F(\xi-s,\eta,\tau;q,M,\tilde M)\,d\eta,
\end{equation}
\begin{equation}\label{calM_3-1}
{\mathcal M}_3(x,t)=-\int_0^t M_2(\tau)\,d\tau\int_0^{t-\tau}M_1(s)\,ds\int_{\frac{s+\tau-t}{2}+x}^x d\xi\int_\tau^{2(\xi-x)+t-s}
F(\xi-s,\eta,\tau;q,M,\tilde M)\,d\eta.
\end{equation}
Taking (\ref{3.6})--(\ref{calQ_3}) and (\ref{calM_1-1})--(\ref{calM_3-1}) into account, we arrive at (\ref{2eq3.4}). $\hfill\Box$

\medskip
Denote $P_q(x,t;M):=P(x,t;q,M)-P(x,t;\tilde q,M)$ and
$$
\|f\|_\delta:=\|f\|_{L_2(0,\delta)}, \quad B_{\delta,r}:=\{f\in L_2(0,\delta):\|f\|_\delta\le r\}.
$$

The following lemma gives estimates for the functions $P(x,t;q,M),$ $P_q(x,t)$ and the function $F(x,t,\tau;q,M,\tilde M),$ determined in
Lemma~6, on their domains of definition as well as an estimate for the function $P_M(x,t),$ determined by (\ref{3.4-1}), on the trapezium
$D_\delta.$

\medskip
{\bf Lemma 7. }{\it For any fixed $r>0,$ the following estimates hold:
\begin{equation}\label{est}
|P(x,t;q,M)|\le C_r, \; |P_q(x,t;M)|\le C_r\|\hat q\|, \; |F(x,t,\tau;q,M,\tilde M)|\le C_r, \; 0\le\tau\le t\le x\le\pi,
\end{equation}
as soon as the functions $q(x),$ $\tilde q(x),$ $(\pi-x)M(x)$ and $(\pi-x)\tilde M(x)$ (which involved) belong to the ball $B_{\pi,r}.$
Moreover, for any $r>0$ and $\delta\in(0,\pi)$ the estimate
\begin{equation}\label{loc}
|P_M(x,t)|\le C_r\sqrt{\delta}\|\hat M\|_\delta, \quad (x,t)\in D_{\delta},
\end{equation}
is fulfilled as soon as $q(x),\,M(x),\,\tilde M(x)\in B_{\delta,r}.$}

\medskip
{\it Proof.} According to Lemmas~2 and~6, respectively, each of the functions $P(x,t;q,M)$ and $F(x,t,\tau;q,M,\tilde M)$ is a solution of
an integral equation having one and the same form (\ref{int_eq}) but with different settings of the parameter $\tau$ and of the free term
$F_0(x,t,\tau).$ Specifically, for $P(x,t;q,M)$ we have $\tau=0$ and $F_0(x,t,0)$ is determined by formula (\ref{ft}), while for
$F(x,t,\tau;q,M,\tilde M)$ the free term is determined by formula (\ref{3.5-1}). Thus, the first and the third estimates in (\ref{est}) are
direct corollaries from (\ref{int_eq_est}), (\ref{FC}) as well as (\ref{ft}) and (\ref{3.5-1}), respectively.

The second estimate in (\ref{est}) can be obtained analogously. Indeed, termwise subtracting equation (\ref{int_eq1}) for the function
$P(x,t;\tilde q,M),$ i.e. after substituting $\tilde q(x)$ instead of $q(x),$ from the unchanged equation (\ref{int_eq1}), one can easily
see that the function $P_q(x,t;M)$ is a solution of equation (\ref{int_eq}) with $\tau=0$ and with the free term
$$
F_0(x,t,0)=\frac12\Big(\int_\frac{t}2^{x-\frac{t}2}\hat q(s)\,ds  +\int_t^x\hat q(s)\,ds\int_0^tP(s,\xi;\tilde q,M)\,d\xi
\qquad\qquad\qquad\qquad\qquad
$$
$$
\;\quad\qquad\qquad +\int_\frac{t}2^t\hat q(s)\,ds\int_0^{2s-t}P(s,\xi;\tilde q,M)\,d\xi -\int_{x-\frac{t}2}^x\hat
q(s)\,ds\int_0^{2(s-x)+t}P(s,\xi;\tilde q,M)\,d\xi\Big),
$$
for which the estimate $|F_0(x,t,0)|\le C_r\|\hat q\|$ holds as soon as $\tilde q(x),\,(\pi-x)M(x)\in B_{\pi,r}.$

It remains to note that estimates (\ref{loc}) follow from (\ref{3.4-1}) along with (\ref{est}) and Observation~1. $\hfill\Box$
\\

{\large\bf 4. Global solvability and uniform full stability of the main equation}
\\

In this section, we give the proof of Theorem~4, stating global solvability of the main nonlinear equation (\ref{ma_eq}) and giving an
estimate for its solution. Then we use the latter for proving the following theorem giving uniform full stability of the main equation.

\medskip
{\bf Theorem 5. }{\it Fix an arbitrary $r>0$ and let the functions $v(x),\,\tilde v(x),\,q(x),\,\tilde q(x)\in B_{\pi,r}$ be given and
satisfy (\ref{means}) along with the analogous condition
\begin{equation}\label{means_ti}
\int_0^\pi \tilde v(x)\,dx =\frac12\int_0^\pi \tilde q(x)\,dx.
\end{equation}
Then the following estimate holds
\begin{equation}\label{est1}
\|\hat M\|_{2,\pi}\le C_r(\|\hat v\|+\|\hat q\|),
\end{equation}
where $M(x)$ is the solution of equation (\ref{ma_eq}), while $\tilde M(x)$ is the one of the equation
\begin{equation}\label{ma_eq_ti}
-\tilde v(\pi-x)=R(\pi,x;\tilde q,\tilde M), \quad 0<x<\pi.
\end{equation}
}

\medskip
Before proceeding directly to the proof of Theorems~4 and~5, we carry out some auxiliary calculations. Consider the function
$F(x,t,\tau;q,M,\tilde M)$ determined in Lemma~6. According to (\ref{int_eq}) and (\ref{3.5-1}), we have
\begin{equation}\label{eq2.2.13}
F(x,t,t;q,M,\tilde M)=\frac{x-t}{2},
\end{equation}
$$
\Phi(x,t;q,M,\tilde M):=\frac{\partial}{\partial x}F(\pi,x,t;q,M,\tilde M)=-\frac12+\frac12\Big(\int_x^\pi P(s-t,x-t;q,\tilde M)\,ds
\qquad\qquad\qquad\qquad\quad
$$
$$
\qquad\qquad\qquad-\int_\frac{x+t}{2}^x P(s-t,2s-x-t;q,\tilde M)\,ds -\int_{\frac{t-x}{2}+\pi}^\pi P(s-t,2(s-\pi)+x-t;q,\tilde M)\,ds\Big)
$$
$$
\quad\qquad\qquad\qquad+\frac12\Big(\int_x^\pi q(s) F(s,x,t;q,M,\tilde M)\,ds -\int_\frac{x+t}{2}^x q(s) F(s,2s-x,t;q,M,\tilde M)\,ds
$$
$$
-\int_{\frac{t-x}{2}+\pi}^\pi q(s) F(s,2(s-\pi)+x,t;q,M,\tilde M)\,ds +\int_0^{x-t} M(s)\,ds\int_x^\pi F(\xi-s,x-s,t;q,M,\tilde M)\,d\xi
$$
$$
-\int_0^{x-t} M(s)\,ds\int_\frac{x+t+s}{2}^x F(\xi-s,2\xi-x-s,t;q,M,\tilde M)\,d\xi
$$
\begin{equation}\label{eq2.2.15}
\qquad\qquad\qquad\qquad\qquad\qquad-\int_0^{x-t} M(s)\,ds\int_{\frac{s+t-x}{2}+\pi}^\pi F(\xi-s,2(\xi-\pi)+x-s,t;q,M,\tilde M)\,d\xi\Big).
\end{equation}
By virtue of the first and the third estimates in (\ref{est}) along with Lemma~2.1. in \cite{But06_MZ}, we have
\begin{equation}\label{est4}
|\Psi(x,t;q,M,\tilde M)|\le f(t), \quad \|f\|\le C_r, \quad \Psi(x,t;q,M,\tilde M):=\frac{2\Phi(x,t;q,M,\tilde M)+1}{\pi-t}.
\end{equation}
as soon as $\|q\|\le r,$ $\|M\|_{2,\pi}\le r$ and $\|\tilde M\|_{2,\pi}\le r.$

\medskip
{\bf Lemma 8. }{\it Fix $r>0$ and choose arbitrary functions $M(x),\,\tilde M(x)\in L_{2,\pi}$ and $\varphi(x)\in L_2(0,\pi),$ satisfying
the conditions $\|M\|_{2,\pi}\le r,$ $\|\tilde M\|_{2,\pi}\le r$ and
\begin{equation}\label{means-1}
\int_0^\pi\varphi(x)\,dx=0.
\end{equation}
Then the solution $h(x)$ of the linear integral equation
\begin{equation}\label{eq2.2.16-0}
\varphi(x)=h(x)-\int_0^x \frac{h(t)}{\pi-t}\,dt +\int_0^x \Psi(x,t;q,M,\tilde M)h(t)\,dt, \quad 0<x<\pi,
\end{equation}
belongs to $L_2(0,\pi)$ and obeys the estimate $\|h\|\le C_r\|\varphi\|.$ }

\medskip
{\it Proof.} Using for an arbitrary function $h(x)\in L_2(0,\pi)$ the mutually inverse transformations
\begin{equation}\label{eq2.2.17-1}
z(x)=h(x)-\int_0^x \frac{h(t)}{\pi-t}\,dt, \quad h(x)=z(x)+\frac1{\pi-x}\int_0^x z(t)\,dt, \quad 0<x<\pi,
\end{equation}
we obtain
\begin{equation}\label{eq2.2.16-3}
 \int_0^x \Psi(x,t;q,M,\tilde M)h(t)\,dt=\int_0^x \Theta(x,t;q,M,\tilde M)z(t)\,dt.
\end{equation}
where
\begin{equation}\label{eq2.2.16-2}
\Theta(x,t;q,M,\tilde M)=\Psi(x,t;q,M,\tilde M)+\int_t^x \frac{\Psi(x,\tau;q,M,\tilde M)}{\pi-\tau}\,d\tau.
\end{equation}
Thus, equation (\ref{eq2.2.16-0}) can be transformed to the equation
\begin{equation}\label{eq2.2.16-1}
\varphi(x)=z(x) +\int_0^x \Theta(x,t;q,M,\tilde M)z(t)\,dt,
\end{equation}
Moreover, by virtue of (\ref{est4}) and (\ref{eq2.2.16-2}) along with Lemma~2.1 in \cite{But06_MZ}, we get
$$
|\Theta(x,t;q,M,\tilde M)|\le f(t)+g(x), \quad \|g\|\le2\|f\|\le C_r.
$$
Hence, in particular, $\|\Theta(\,\cdot\,\,\cdot\,;q,M,\tilde M)\|_{L_2((0,\pi)^2)}\le C_r.$ Thus, equation (\ref{eq2.2.16-1}) has a unique
square integrable solution  $z(x),$ which, by virtue of Lemma~1 in \cite{ButMal18}, obeys the estimate $\|z\|\le C_r\|\varphi\|.$

Let us show that
\begin{equation}\label{means-2}
\int_0^\pi dx\int_0^x \Theta(x,t;q,M,\tilde M)z(t)\,dt=0.
\end{equation}
Indeed, according to the second relation in (\ref{eq2.2.17-1}) along with square integrability of $z(x),$ we have $(\pi-x)^\theta h(x)\in
L_2(0,\pi)$ for any $\theta>1/2.$ Further, by virtue of (\ref{eq2.2.13})--(\ref{est4}) and (\ref{eq2.2.16-3}), we get
$$
\theta(\zeta):=\int_0^\zeta dx\int_0^x \Theta(x,t;q,M,\tilde M)z(t)\,dt=\int_0^\zeta dx\int_0^x \Big(2\frac{\partial}{\partial x}
F(\pi,x,\tau;q,M,\tilde M)+ 1\Big)\frac{h(\tau)}{\pi-\tau}\,d\tau \qquad\qquad\qquad
$$
$$
=\int_0^\zeta \frac{h(\tau)}{\pi-\tau}\,d\tau\int_\tau^\zeta \Big(2\frac{\partial}{\partial x} F(\pi,x,\tau;q,M,\tilde
M)+1\Big)\,dx=\int_0^\zeta \frac{h(\tau)}{\pi-\tau}\Big(2F(\pi,\zeta,\tau;q,M,\tilde M)-\pi+\zeta\Big)\,d\tau.
$$
Thus, according to (\ref{int_eq}) and (\ref{3.5-1}), we arrive at the relation
\begin{equation}\label{2eq5.8}
\theta(\zeta)=\sum_{k=1}^6\int_0^\zeta \frac{h(\tau)}{\pi-\tau}F_k(\zeta,\tau)\,d\tau,
\end{equation}
where
$$
F_1(\zeta,\tau):=\int_\zeta^\pi ds\int_0^{\zeta-\tau} P(s-\tau,\xi;q,\tilde M)\,d\xi,
$$
$$
F_2(\zeta,\tau):=\int_\frac{\zeta+\tau}{2}^\zeta ds\int_0^{2s-\zeta-\tau}P(s-\tau,\xi;q,\tilde M)
\,d\xi-\int_{\frac{\tau-\zeta}{2}+\pi}^\pi ds\int_0^{2(s-\pi)+\zeta-\tau} P(s-\tau,\xi;q,\tilde M)\,d\xi,
$$
$$
F_3(\zeta,\tau):=\int_\zeta^\pi q(s)\,ds\int_\tau^\zeta F(s,\xi,\tau;q,M,\tilde M)\,d\xi,
$$
$$
F_4(\zeta,\tau):=\int_\frac{\zeta+\tau}{2}^\zeta q(s)\,ds\int_\tau^{2s-\zeta} F(s,\xi,\tau;q,M,\tilde
M)\,d\xi-\int_{\frac{\tau-\zeta}{2}+\pi}^\pi q(s)\,ds\int_\tau^{2(s-\pi)+\zeta} F(s,\xi,\tau;q,M,\tilde M)\,d\xi,
$$
$$
F_5(\zeta,\tau):=\int_0^{\zeta-\tau} M(s)\,ds\int_\zeta^\pi d\xi\int_\tau^{\zeta-s} F(\xi-s,\eta,\tau;q,M,\tilde M)\,d\eta,
$$
$$
F_6(\zeta,\tau):=\int_0^{\zeta-\tau} M(s)\,ds\int_\frac{\zeta+\tau+s}{2}^\zeta d\xi\int_\tau^{2\xi-\zeta-s}F(\xi-s,\eta,\tau;q,M,\tilde
M)\,d\eta\qquad\qquad\qquad\qquad
$$
$$
\qquad\qquad\qquad\qquad-\int_0^{\zeta-\tau} M(s)\,ds\int_{\frac{s+\tau-\zeta}{2}+\pi}^\pi d\xi\int_\tau^{2(\xi-\pi)+\zeta-s}
F(\xi-s,\eta,\tau;q,M,\tilde M)\,d\eta.
$$
By virtue of the first estimate in (\ref{est}), we have
$$
|F_1(\zeta,\tau)|\le \int_\zeta^\pi ds\int_0^{\zeta-\tau} |P(s-\tau,\xi;q,\tilde M)|\,d\xi\le C_r(\pi-\zeta)(\zeta-\tau).
$$
Moreover, continuing the function $P(x,t;q,\tilde M)$ by zero outside the triangle $D_\pi$ and transforming the limits of integration, we
get
$$
F_2(\zeta,\tau)=\int_\frac{\zeta+\tau}{2}^{\frac{\tau-\zeta}{2}+ \pi} ds\int_0^{2s-\zeta-\tau}P(s-\tau,\xi;q,\tilde M) \,d\xi
+\int_{\frac{\tau-\zeta}{2}+\pi}^\zeta ds\int_{2(s-\pi)+\zeta-\tau}^{2s-\zeta-\tau}P(s-\tau,\xi;q,\tilde M)\,d\xi
$$
$$
\;\quad\qquad\qquad\qquad\qquad\qquad\qquad\qquad\qquad\qquad\qquad -\int_\zeta^\pi ds\int_0^{2(s-\pi)+\zeta-\tau} P(s-\tau,\xi;q,\tilde
M)\,d\xi,
$$
where
$$
\Big|\int_\frac{\zeta+\tau}{2}^{\frac{\tau-\zeta}{2}+ \pi} ds\int_0^{2s-\zeta-\tau}P(s-\tau,\xi;q,\tilde M)\,d\xi\Big|\le C_r
\int_\frac{\zeta+\tau}{2}^{\frac{\tau-\zeta}{2}+ \pi} ds\int_0^{2s-\zeta-\tau}d\xi \le C_r(\pi-\zeta)^2,
$$
$$
\Big|\int_{\frac{\tau-\zeta}{2}+\pi}^\zeta ds\int_{2(s-\pi)+\zeta-\tau}^{2s-\zeta-\tau}P(s-\tau,\xi;q,\tilde M)\,d\xi\Big|\le
C_r\Big(\frac{\zeta-\tau}2+\pi-\zeta\Big)(\pi-\zeta),
$$
$$
\Big|\int_\zeta^\pi ds\int_0^{2(s-\pi)+\zeta-\tau} P(s-\tau,\xi;q,\tilde M)\,d\xi\Big|\le C_r(\pi-\zeta)\Big(2(\pi-\zeta)+\zeta-\tau\Big).
$$
Thus, we arrive at the estimate
$$
|F_2(\zeta,\tau)|\le C_r(\pi-\zeta)(\pi-\tau).
$$
Analogously, using the third estimate in  (\ref{est}), we get the following estimates for other $F_k(\zeta,\tau)$'s:
$$
|F_3(\zeta,\tau)|\le \int_\zeta^\pi |q(s)|\,ds\int_\tau^\zeta |F(s,\xi,\tau;q,M,\tilde M)|\,d\xi \le C_r\sqrt{\pi-\zeta}(\zeta-\tau),
$$
$$
|F_4(\zeta,\tau)|\le\int_\frac{\zeta+\tau}{2}^{\frac{\tau-\zeta}{2}+ \pi}|q(s)|\, ds\int_\tau^{2s-\zeta}|F(s,\xi,\tau;q,M,\tilde M)|\,d\xi
\qquad\qquad\qquad\quad\qquad\qquad\qquad\qquad\qquad\qquad
$$
$$
+\Big|\int_{\frac{\tau-\zeta}{2}+\pi}^\zeta q(s)\, ds\int_{2(s-\pi)+\zeta}^{2s-\zeta}F(s,\xi,\tau;q,M,\tilde M)\,d\xi\Big|
\qquad\qquad\qquad
$$
$$
\qquad\qquad\qquad\qquad\qquad +\Big|\int_\zeta^\pi q(s)\, ds\int_\tau^{2(s-\pi)+\zeta} F(s,\xi,\tau;q,M,\tilde M)\,d\xi\Big|\le
C_r\sqrt{\pi-\zeta}(\pi-\tau)
$$
$$
|F_5(\zeta,\tau)|\le\int_0^{\zeta-\tau} |M(s)|\,ds\int_\zeta^\pi d\xi\int_\tau^{\zeta-s} |F(\xi-s,\eta,\tau;q,M,\tilde M)|\,d\eta\le
C_r(\pi-\zeta)\sqrt{\zeta-\tau},
$$
$$
|F_6(\zeta,\tau)|\le\int_0^{\zeta-\tau}|M(s)|\,ds \int_\frac{\zeta+\tau+s}{2}^{\frac{s+\tau-\zeta}{2}+ \pi} d\xi
\int_\tau^{2\xi-\zeta-s}|F(\xi-s,\eta,\tau;q,M,\tilde M)|\,d\eta \qquad\qquad\qquad\qquad\qquad\qquad
$$
$$
+\Big|\int_0^{\zeta-\tau}M(s)\,ds \int_{\frac{s+\tau-\zeta}{2}+\pi}^\zeta d\xi
\int_{2(\xi-\pi)+\zeta-s}^{2\xi-\zeta-s}F(\xi-s,\eta,\tau;q,M,\tilde M)\,d\eta\Big| \qquad
$$
$$
\quad\qquad\quad +\Big|\int_0^{\zeta-\tau}M(s)\,ds \int_\zeta^\pi d\xi\int_\tau^{2(\xi-\pi)+\zeta-s} F(\xi-s,\eta,\tau;q,M,\tilde
M)\,d\eta\Big|\le C_r(\pi-\zeta)\sqrt{\zeta-\tau}.
$$
Thus, the following common estimate holds:
$$
|F_k(\zeta,\tau)|\le C_r\sqrt{\pi-\zeta}\,(\pi-\tau), \quad k=\overline{1,6},
$$
which along with (\ref{2eq5.8}) and $(\pi-x)^\theta h(x)\in L_2(0,\pi),$ $\theta\in(1/2,3/4),$ give
$$
|\theta(\zeta)| \le C_r\sqrt{\pi-\zeta} \int_0^\zeta |h(\tau)|\,d\tau \le C_r(\pi-\zeta)^\frac14 \int_0^\zeta
\frac{(\pi-\tau)^\theta|h(\tau)|}{(\pi-\tau)^{\theta-\frac14}}\,d\tau\to0, \quad \zeta\to\pi,
$$
which implies (\ref{means-2}). Further, according to (\ref{means-1}), (\ref{eq2.2.16-1}) and (\ref{means-2}), we have
$$
\int_0^\pi z(x)\,dx=0,
$$
which along with the second equality in (\ref{eq2.2.17-1}) give
$$
h(x)=z(x)-\frac1{\pi-x}\int_x^\pi z(t)\,dt, \quad 0<x<\pi,
$$
By virtue of Lemma~2.1 in \cite{But06_MZ}, we get $\|h\|\le3\|z\|\le C_r\|\varphi\|,$ which finishes the proof. $\hfill\Box$

\medskip
Now we are in position to give the proof of Theorems~4 and~5.

\medskip
{\it Proof of Theorem~4.} By virtue of (\ref{int_eq1}) and (\ref{K_12}), we have
$$
R(x,t;q,M)=\frac12\Big(-\frac12\Big(q\Big(x-\frac{t}{2}\Big)+q\Big(\frac{t}{2}\Big)\Big)+(x-t)M(t)-\int_0^t M(s)\,ds
\qquad\qquad\qquad\qquad
$$
$$
\quad\qquad\qquad+\int_t^x q(s)P(s,t)\,ds-\int_\frac{t}{2}^t q(s)P(s,2s-t)\,ds-\int_{x-\frac{t}{2}}^x q(s)P(s,2(s-x)+t)\,ds
$$
$$
\;\quad\qquad\qquad+\int_0^t M(s)\,ds\int_t^x P(\xi-s,t-s)\,d\xi - \int_0^t M(s)\,ds\int_\frac{t+s}{2}^t P(\xi-s,2\xi-t-s)\,d\xi
$$
\begin{equation}\label{R}
\qquad\qquad\qquad\qquad\qquad\qquad\qquad\qquad\qquad-\int_0^t M(s)\,ds\int_{\frac{s-t}{2}+x}^x P(\xi-s,2(\xi-x)+t-s)\,d\xi\Big),
\end{equation}
where $P(x,t)=P(x,t;q,M).$ Thus, equation (\ref{ma_eq}) can be represented in the form
\begin{equation}\label{ma_eq1}
M(x)=g(x)+{\cal D}_qM(x), \quad 0<x<\pi,
\end{equation}
where
\begin{equation}\label{ft1}
g(x)=\frac{1}{\pi-x}\Big(\frac12\Big(q\Big(\pi-\frac{x}{2}\Big)+q\Big(\frac{x}{2}\Big)\Big)-2v(\pi-x)\Big),
\end{equation}
$$
{\cal D}_qM(x)=\frac{1}{\pi-x}\Big(\int_0^x M(s)\,ds - \int_x^\pi q(s)P(s,x;q,M)\,ds+\int_\frac{x}{2}^x q(s)P(s,2s-x;q,M)\,ds\quad
$$
$$
\qquad\quad+\int_{\pi-\frac{x}{2}}^\pi q(s)P(s,2(s-\pi)+x;q,M)\,ds- \int_0^x M(s) \Big(\int_x^\pi P(\xi-s,x-s;q,M)\,d\xi
$$
\begin{equation}\label{DqM}
\;\;\qquad\qquad\quad -\int_\frac{x+s}{2}^x P(\xi-s,2\xi-x-s;q,M)\,d\xi -\int_{\frac{s-x}{2}+\pi}^\pi
P(\xi-s,2(\xi-\pi)+x-s;q,M)\,d\xi\Big)ds\Big).
\end{equation}
The latter implies
$$
{\cal D}_qM(x)-{\cal D}_q\tilde M(x) =\frac{1}{\pi-x}\Big(\int_0^x \hat M(s)\,ds - \int_x^\pi q(s)P_M(s,x)\,ds+\int_\frac{x}{2}^x
q(s)P_M(s,2s-x)\,ds \qquad\;
$$
$$
\qquad\qquad\qquad\qquad+\int_{\pi-\frac{x}{2}}^\pi q(s)P_M(s,2(s-\pi)+x)\,ds- \int_0^x \hat M(s) \Big(\int_x^\pi P(\xi-s,x-s;q,M)\,d\xi
$$
$$
\qquad\qquad\qquad\qquad\quad -\int_\frac{x+s}{2}^x P(\xi-s,2\xi-x-s;q,M)\,d\xi -\int_{\frac{s-x}2+\pi}^\pi
P(\xi-s,2(\xi-\pi)+x-s;q,M)\,d\xi\Big)ds\Big)
$$
$$
\quad\qquad\qquad- \int_0^x \tilde M(s) \Big(\int_x^\pi P_M(\xi-s,x-s)\,d\xi -\int_\frac{x+s}2^x P_M(\xi-s,2\xi-x-s)\,d\xi
$$
\begin{equation}\label{DqM_hat}
\quad\quad\quad\qquad\qquad\qquad\qquad\qquad\qquad\qquad\qquad\qquad-\int_{\frac{s-x}2+\pi}^\pi
P_M(\xi-s,2(\xi-\pi)+x-s)\,d\xi\Big)ds\Big).
\end{equation}
Let $\|q\|\le r$ and $\|v\|\le r.$ Then, according to (\ref{ft1}), for any $\delta\in(0,\pi]$ we have $\|g\|_\delta \le\|g\|\le\alpha_r,$
where $\alpha_r$ depends only on $r.$ Further, by virtue of (\ref{DqM}) and (\ref{DqM_hat}) along with Lemma~7, we arrive at the estimates
$$
 \quad \|{\cal D}_qM\|_{\delta_0}\le\beta_r\sqrt{\delta_0}, \quad \|{\cal D}_qM-{\cal D}_q\tilde M\|_{\delta_0}\le\gamma_r\delta_0\|\hat M\|_{\delta_0}
$$
for any $\delta_0\in(0,\pi)$ as soon as $\|M\|_{\delta_0}\le2\alpha_r$ and $\|\tilde M\|_{\delta_0}\le2\alpha_r,$ where $\beta_r$ and
$\gamma_r$ depend only on~$r.$ Choose $\delta_0=\delta_0(r)\in(0,\pi)$ so that $\beta_r\sqrt{\delta_0}\le\alpha_r$ and
$\gamma_r\delta_0<1.$ Then for any $\delta\in(0,\delta_0],$ the operator
$$
{\cal A}_qM(x):=g(x)+{\cal D}_qM(x), \quad 0<x<\delta,
$$
maps the ball $B_{\delta,2\alpha_r}$ into itself and is a contraction mapping in $B_{\delta,2\alpha_r}.$ Thus, according to the contraction
mapping theorem, equation (\ref{ma_eq1}) has a unique solution $M(x)$ on the interval $(0,\delta)$ that belongs to the ball
$B_{\delta,2\alpha_r}.$

Let us have already found a solution $M_1(x)$ of the main equation (\ref{ma_eq}) on the interval $(0,\delta)$ for a certain
$\delta\in(0,\pi/2]$ and let $M_1(x)\in B_{\delta,C_r}.$ Then we continue the function $M_1(x)$ to $(\delta,2\delta)$ by zero and seek a
solution on $(0,2\delta)$ in the form $M_1(x)+M_2(x),$ where $M_2(x)=0$ on $(0,\delta).$ By virtue of (\ref{ma_eq}), (\ref{K_12}),
(\ref{2eq3.4}), (\ref{eq2.2.13}) and (\ref{eq2.2.15}), we arrive at the following linear integral equation with respect to the function
$M_2(x)$ on $(\delta,2\delta):$
\begin{equation}\label{eq2.2.14}
\varphi(x)=\frac{\pi-x}{2}M_2(x)+\int_\delta^x \Phi(x,t;q,M_1,M_1)M_2(t)\,dt, \quad \delta<x<2\delta,
\end{equation}
where
\begin{equation}\label{ft2}
\varphi(x)=-v(\pi-x)-R(\pi,x;q,M_1), \quad 0<x<2\delta.
\end{equation}
Using the first estimate in (\ref{est}) and (\ref{R}), (\ref{ft2}) as well as (\ref{est4}), we get
\begin{equation}\label{est2}
\|\varphi\|_{L_2(\delta,2\delta)}\le C_r, \quad \|\Phi(\,\cdot\,,\,\cdot\,;q,M_1,M_1)\|_{L_2((\delta,2\delta)^2)}\le C_r.
\end{equation}
Hence, equation (\ref{eq2.2.14}) has a unique solution $M_2(x).$ Obviously, the function $M(x):=M_1(x)+M_2(x)$ is the unique solution of
equation (\ref{ma_eq}) on $(0,2\delta)$ that coincides with $M_1(x)$ on $(0,\delta).$ Moreover, if $2\delta<\pi,$ then $M_2(x)\in
L_2(\delta,2\delta)$ and, by virtue of Lemma~1 in \cite{ButMal18}, we get $M(x)\in B_{2\delta,C_r}.$

Continuing the process, we obtain a solution $M(x)$ of the main equation (\ref{ma_eq}) on the entire interval $(0,\pi)$ belonging to
$L_2(0,T)$ for each $T\in(0,\delta)$ and, in particular, $\|M\|_{\pi/2}\le C_r.$ This solution is unique. Indeed, let there exists another
locally square integrable solution $\tilde M(x).$ Then for sufficiently small $\delta\in(0,\delta_0(r)]$ both the functions $M(x)$ and
$\tilde M(x)$ on $(0,\delta)$ belong to the ball $B_{\delta,2\alpha_r}$ and, according to the first part of the proof, they coincide a.e.
on $(0,\delta).$ Then, by virtue of uniqueness of the continuation of solution, they coincide a.e. on $(0,\pi).$

It remains to prove that $M(x)\in L_{2,\pi}$ and $\|M\|_{2,\pi}\le C_r.$ For this purpose, we rewrite the linear equation (\ref{eq2.2.14})
for $\delta=\pi/2$ in the form
$$
2\varphi(x)=h(x)-\int_0^x \frac{h(t)}{\pi-t}\,dt +\int_0^x \Psi(x,t;q,M_1,M_1)h(t)\,dt, \quad 0<x<\pi,
$$
where $h(x)=(\pi-x)M_2(x)$ and the function $\Psi(x,t;q,M,\tilde M)$ is determined in (\ref{est4}). Note that we painlessly extended the
interval $(\pi/2,\pi)$ till $(0,\pi),$ because $h(x)=0$ on $(0,\pi/2).$ By virtue of (\ref{means}), the function $\varphi(x)$ obeys
(\ref{means-1}). Indeed, subsequently using (\ref{ft2}), (\ref{means}), (\ref{K_12}) and (\ref{Kx0}), we calculate
$$
\int_0^\pi\varphi(x)\,dx= -\int_0^\pi v(x)\,dx-\int_0^\pi R(\pi,x;q,M_1)\,dx \qquad\qquad\qquad
$$
$$
\;\;\qquad\qquad\qquad =-\frac12\int_0^\pi q(x)\,dx -P(\pi,\pi;q,M_1) +P(\pi,0;q,M_1)=0.
$$
By virtue of Lemma~8 and (\ref{est2}), we get $h(x)\in L_2(0,\pi)$ and $\|h\|\le C_r,$ which finishes the proof. $\hfill\Box$

\medskip
{\it Proof of Theorem 5.} First, let us prove the theorem in the particular case when $\tilde q(x)=q(x).$ Termwise subtracting
(\ref{ma_eq_ti}) from (\ref{ma_eq}) and using (\ref{K_12}) along with Lemma~6, we get
\begin{equation}\label{eq2.2.16-00}
\frac12\varphi(x) =R(\pi,x;q,M) -R(\pi,x;q,\tilde M) =\frac{d}{dx}\int_0^x F(\pi,x,t;q,M,\tilde M)\hat M(t)\,dt,
\end{equation}
where $\varphi(x)=-2\hat v(\pi-x),$ which, according to (\ref{means}), (\ref{means_ti}) and our assumption that $\tilde q(x)=q(x),$
satisfies (\ref{means-1}). By virtue of (\ref{eq2.2.13})--(\ref{est4}), relation (\ref{eq2.2.16-00}) can be rewritten in the form
(\ref{eq2.2.16-0}) with $h(x)=(\pi-x)\hat M(x).$  Using Lemma~8 along with the second part of Theorem~4, we obtain $\|h\|\le
C_r\|\varphi\|,$ which coincide with estimate (\ref{est1}) for $\hat q(x)=0.$

In the general case, when $\tilde q(x)$ may differ from $q(x),$ we rewrite (\ref{ma_eq_ti}) in the form
\begin{equation}\label{ma_eq_ti-1}
v_1(x)=R(\pi,x;q,\tilde M),
\end{equation}
where $v_1(x)=-\tilde v(\pi-x) +R(\pi,x;q,\tilde M) -R(\pi,x;\tilde q,\tilde M).$ Subtracting (\ref{ma_eq_ti-1}) from (\ref{ma_eq}), we get
relation (\ref{eq2.2.16-00}) with $\varphi(x)=-2\hat v(\pi-x) -2R(\pi,x;q,\tilde M) +2R(\pi,x;\tilde q,\tilde M).$ By virtue of
(\ref{Kx0}), (\ref{K_12}), (\ref{means}) and (\ref{means_ti}), we have
$$
\frac12\int_0^\pi\varphi(x)\,dx= -\int_0^\pi \hat v(x)\,dx -\int_0^\pi R(\pi,x;q,\tilde M)\,dx +\int_0^\pi R(\pi,x;\tilde q,\tilde M)\,dx
\qquad\qquad\qquad\qquad\quad
$$
$$
\qquad\qquad\qquad =-\frac12\int_0^\pi \hat q(x)\,dx -P(\pi,\pi;q,\tilde M) +P(\pi,0;q,\tilde M) +P(\pi,\pi;\tilde q,\tilde M)
-P(\pi,0;\tilde q,\tilde M)
$$
$$
\qquad\qquad\qquad\qquad\qquad\qquad\qquad\qquad\qquad=- \frac12\int_0^\pi \hat q(x)\,dx +\frac12\int_0^\pi q(x)\,dx -\frac12\int_0^\pi
\tilde q(x)\,dx =0.
$$
Thus, according to the first part of the proof, we obtain the estimate $\|\hat M\|_{2,\pi}\le C_r\|\varphi\|,$ where $\|\varphi\| \le
C_r(\|\hat v\| +\|p\|)$ and $p(x)=R(\pi,x;q,\tilde M) -R(\pi,x;\tilde q,\tilde M).$ Thus, it remains to show that
\begin{equation}\label{est5}
\|p\|\le C_r\|\hat q\|.
\end{equation}
By virtue of (\ref{R}), we have
$$
p(x)=\frac12\Big(-\frac12\Big(\hat q\Big(\pi-\frac{x}{2}\Big)+\hat q\Big(\frac{x}{2}\Big)\Big)+\int_x^\pi \hat q(s)P(s,x;q,\tilde
M)\,ds-\int_\frac{x}{2}^x \hat q(s)P(s,2s-x;q,\tilde M)\,ds
$$
$$
-\int_{\pi-\frac{x}{2}}^\pi \hat q(s)P(s,2(s-\pi)+x;q,\tilde M)\,ds+\int_x^\pi \tilde  q(s)P_q(s,x;\tilde M)\,ds-\int_\frac{x}{2}^x \tilde
q(s)P_q(s,2s-x;\tilde M)\,ds
$$
$$
-\int_{\pi-\frac{x}{2}}^\pi \tilde q(s)P_q(s,2(s-\pi)+x;\tilde M)\,ds +\int_0^x \tilde M(s)\,ds\int_x^\pi P_q(\xi-s,x-s;\tilde M)\,d\xi
$$
$$
- \int_0^x \tilde M(s)\,ds\int_\frac{x+s}{2}^x P_q(\xi-s,2\xi-x-s;\tilde M)\,d\xi -\int_0^x \tilde M(s)\,ds\int_{\frac{s-x}{2}+\pi}^\pi
P_q(\xi-s,2(\xi-\pi)+x-s;\tilde M)\,d\xi\Big),
$$
Using the first and the second estimates in (\ref{est}), we obtain (\ref{est5}). $\hfill\Box$
\\

{\bf \large 5. Uniform stability of recovering the characteristic function kernel}\\

Let $\lambda_n,$ $n\ge1,$ be all zeros with account of multiplicity of an entire function $\Delta(\lambda)$ of the form
\begin{equation}\label{char1-1}
\Delta(\lambda)=\frac{\sin\rho\pi}{\rho} +\int_0^\pi v(x)\frac{\cos\rho x}{\rho^2}\,dx, \quad \rho^2=\lambda, \quad v(x)\in L_2(0,\pi),
\quad \int_0^\pi v(x)\,dx=0.
\end{equation}
As was mentioned in Section~2, the numbers $\lambda_n$ have asymptotics (\ref{asy}) with $\omega=0,$ and the function $\Delta(\lambda)$ is
determined by them uniquely by formula (\ref{char2}). Conversely, the function $\Delta(\lambda),$ constructed by formula (\ref{char2}) from
any sequence $\{\lambda_n\}_{n\ge1}$ of complex numbers of the form (\ref{asy}) with $\omega=0,$ has the form (\ref{char1-1}) with some
square integrable function $v(x)$ possessing zero mean value on $(0,\pi).$ Consider another sequence
$$
\tilde\lambda_n=n^2 +\tilde\varkappa_n, \quad n\ge1, \quad \{\tilde\varkappa_n\}\in l_2,
$$
along with the corresponding functions
\begin{equation}\label{tilde_Delta}
\tilde \Delta(\lambda)=\pi\prod_{n=1}^\infty\frac{\tilde\lambda_n-\lambda}{n^2} =\frac{\sin\rho\pi}\rho +\int_0^\pi \tilde
v(x)\frac{\cos\rho x}{\rho^2}\,dx, \quad \tilde v(x)\in L_2(0,\pi), \quad \int_0^\pi \tilde v(x)\,dx=0.
\end{equation}
In this section, we prove the following theorem, which gives uniform stability of recovering the kernel $v(x)$ from zeros
$\{\lambda_n\}_{n\ge1}$ of the function $\Delta(\lambda).$

\medskip
{\bf Theorem 6. }{\it For any $r>0,$ the estimate
\begin{equation}\label{est_w}
\|\hat v\|\le C_r\|\{\hat\varkappa_n\}\|_{l_2}
\end{equation}
is fulfilled as soon as $\|\{\varkappa_n\}\|_{l_2}\le r$ and $\|\{\tilde\varkappa_n\}\|_{l_2}\le r.$}

\medskip
Before proceeding directly to the proof of Theorem~6, we provide several auxiliary assertions. First of all, we prove (\ref{est_w}) in the
particular case when $\tilde\lambda_n=n^2,$ $n\ge1,$ i.e. $\tilde v(x)=0.$

\medskip
{\bf Proposition 1. }{\it For any $r>0,$ the estimate $\|v\|\le C_r\|\{\varkappa_n\}\|_{l_2}$ holds as soon as
$\|\{\varkappa_n\}\|_{l_2}\le r.$}

\medskip
{\it Proof.} According to (\ref{char2}), (\ref{char1-1}) and Parseval's equality, we calculate
$$
\|v\|=\sqrt{\frac2\pi\sum_{k=1}^\infty k^4|\Delta(k^2)|^2}, \quad k^2\Delta(k^2) =\pi k^2\prod_{j=1}^\infty \frac{\lambda_j-k^2}{j^2}= a_k
b_k\varkappa_k, \quad a_k=\prod_{{j\ne k}\atop{j=1}}^\infty \frac{\lambda_j-k^2}{j^2-k^2},
$$
$$
b_k=\pi\prod_{{j\ne k}\atop{j=1}}^\infty \frac{j^2-k^2}{j^2} =k^2\lim_{\rho\to k}\frac{\sin\rho\pi}{\rho(k+\rho)(k-\rho)}
=-\frac\pi2\lim_{\rho\to k}\cos\rho\pi =(-1)^{k+1}\frac\pi2.
$$
Thus, one needs to prove that $|a_k|\le C_r$ uniformly as $\|\{\varkappa_n\}\|_{l_2}\le r.$ For this purpose, we factorize
$a_k=a_{k,1}a_{k,2},$ where $a_{k,2}=1$ for $r\le1$ and
$$
a_{k,1}=\prod_{|j-k|\ge r}\Big(1+\frac{\varkappa_j}{j^2-k^2}\Big), \quad a_{k,2}=\prod_{0<|j-k|<r}\Big(1+\frac{\varkappa_j}{j^2-k^2}\Big),
\quad r>1.
$$
Since
$$
\frac{|\varkappa_j|}{|j^2-k^2|} =\frac{|\varkappa_j|}{|j-k|(j+k)} \le\frac{\|\{\varkappa_n\}\|_{l_2}}{2r} \le\frac12, \quad |j-k|\ge r,
$$
using the Cauchy--Bunyakovsky--Schwarz inequality, we get the estimate
$$
|a_{k,1}|=\Big|\exp\Big(\sum_{|j-k|\ge r}\ln\Big(1+\frac{\varkappa_j}{j^2-k^2}\Big)\Big)\Big|\le \exp\Big(2\sum_{|j-k|\ge
r}\frac{|\varkappa_j|}{|j^2-k^2|}\Big) <\exp\left(2\sqrt{\sum_{j=2}^\infty\frac1{j^2}}\right).
$$
Finally, we have
$$
|a_{k,2}|\le\prod_{0<|j-k|< r}\Big(1+\frac{r}{|j^2-k^2|}\Big) <(1+r)^{2r-2}, \quad r>1,
$$
which finishes the proof. $\hfill\Box$

\medskip
In what follows, without loss of generality, we assume $r\in{\mathbb N}.$ For $k\in{\mathbb N}$ we introduce the sets
$$
\Omega_r(k):=\{j:|j-k|<3r,j\in{\mathbb N}\}, \;\; \Omega_r'(k):=\Omega_r(k)\setminus\{k\}, \;\; \Theta_r(k):=\{j:|j-k|\ge 3r,j\in{\mathbb
N}\}.
$$
Obviously, $\Omega_r(k)\cup\Theta_r(k)\equiv{\mathbb N}.$ We also put $\alpha_{r,k}:=\#\Omega_r(k)=\min\{k,3r\}+3r-1\le 6r-1.$

\medskip
Denote
$$
\sigma_{r,k}(\lambda):=\prod_{j\in\Omega_r(k)}\frac{\lambda_j-\lambda}{j^2}, \quad
\tilde\sigma_{r,k}(\lambda):=\prod_{j\in\Omega_r(k)}\frac{\tilde\lambda_j-\lambda}{j^2}.
$$

\medskip
{\bf Proposition 2. }{\it For any $r\in{\mathbb N},$ the estimates
\begin{equation}\label{theta_est}
|\sigma_{r,k}(k^2)|\le C_r\frac{|\varkappa_k|}{k^{\alpha_{r,k}+1}}, \quad |\sigma_{r,k}(k^2)
-\tilde\sigma_{r,k}(k^2)|\le\frac{C_r}{k^{\alpha_{r,k}+1}}\sum_{j\in\Omega_r(k)}|\hat\varkappa_j|, \quad k\in{\mathbb N},
\end{equation}
are fulfilled as soon as $\|\{\varkappa_n\}\|_{l_2}\le r$ and $\|\{\tilde\varkappa_n\}\|_{l_2}\le r.$}

\medskip
{\it Proof.} We have
\begin{equation}\label{theta_est1}
\sigma_{r,k}(k^2)=\frac{\varkappa_k}{k^2}\prod_{j\in\Omega_r'(k)}\frac{\lambda_j-k^2}{j^2}.
\end{equation}
Since $j\in\Omega_r(k)$ is equivalent to the inequalities $\max\{0,k-3r\}< j< k+3r,$ we have
\begin{equation}\label{theta_est2}
|\lambda_j-k^2| =|j^2-k^2+\varkappa_j| < 3r(2k+3r) +r \le C_rk, \quad j\in\Omega_r(k), \;\; k\in{\mathbb N},
\end{equation}
and
\begin{equation}\label{theta_est3}
\frac1j\le\frac1{\max\{0,k-3r\}+1} \le\frac{3r}{k}, \quad j\in\Omega_r(k), \;\; k\in{\mathbb N}.
\end{equation}
Substituting estimates (\ref{theta_est2}) and (\ref{theta_est3}) into (\ref{theta_est1}), we get
$$
|\sigma_{r,k}(k^2)|\le \frac{|\varkappa_k|}{k^2} \Big(\frac{9r^2C_r}k\Big)^{\alpha_{r,k}-1},
$$
which coincides with the first estimate in (\ref{theta_est}). Further, it is easy to check that
\begin{equation}\label{theta_est4}
\sigma_{r,k}(k^2) -\tilde\sigma_{r,k}(k^2) =\sum_{j\in\Omega_r(k)}\tilde\sigma_{r,k,j}(k^2)\frac{\hat\varkappa_j}{j^2}\sigma_{r,k,j}(k^2),
\end{equation}
where
$$
\sigma_{r,k,j}(\lambda)=\prod_{\nu=j+1}^{k+3r-1}\frac{\lambda_\nu-\lambda}{\nu^2}, \quad
\tilde\sigma_{r,k,j}(\lambda)=\prod_{\nu=\max\{0,k-3r\}+1}^{j-1}\frac{\tilde\lambda_\nu-\lambda}{\nu^2}.
$$
According to (\ref{theta_est2}) and (\ref{theta_est3}), we have
$$
|\sigma_{r,k,j}(k^2)\tilde\sigma_{r,k,j}(k^2)| \le\Big(\frac{9r^2C_r}k\Big)^{\alpha_{r,k}-1},
$$
which along with (\ref{theta_est3}) and (\ref{theta_est4}) give the second estimate in (\ref{theta_est}). $\hfill\Box$

\medskip
Denote
\begin{equation}\label{Delta_k}
\Delta_k(\lambda):=\frac{\Delta(\lambda)}{\sigma_{r,k}(\lambda)} =\pi\prod_{j\in\Theta_r(k)}\frac{\lambda_j-\lambda}{j^2}, \quad
\tilde\Delta_k(\lambda):=\frac{\tilde\Delta(\lambda)}{\tilde\sigma_{r,k}(\lambda)}
=\pi\prod_{j\in\Theta_r(k)}\frac{\tilde\lambda_j-\lambda}{j^2}.
\end{equation}

\medskip
{\bf Proposition 3. }{\it For any $r\in{\mathbb N},$ the estimate
\begin{equation}\label{Delta_k-0}
|\Delta_k(k^2)|\le C_r k^{\alpha_{r,k}-1}, \quad k\in{\mathbb N},
\end{equation}
holds as soon as $\|\{\varkappa_n\}\|_{l_2}\le r.$}

\medskip
{\it Proof.} According to (\ref{asy}) with $\omega=0,$ we have
$$
\lambda_n=\rho_n^2, \quad \rho_n=n+\frac{\varepsilon_n}n, \quad |\varepsilon_n|\le|\varkappa_n|, \quad n\ge1.
$$
Indeed, it is easy to see that $\varepsilon_n=n^2(\sqrt{1+\varkappa_n/n^2}-1),$ where ${\rm Re}\sqrt{\,\cdot\,}\ge0.$ Hence,
$|\varepsilon_n|= |\varkappa_n|/|\sqrt{1+\varkappa_n/n^2} +1| \le |\varkappa_n|/({\rm Re}\sqrt{1+\varkappa_n/n^2} +1) \le|\varkappa_n|.$

Thus, assuming $\|\{\varkappa_n\}\|_{l_2}\le r,$ we get $|\varepsilon_n|\le r.$ Put $\xi_r:= 5r+1.$ Then, in particular, we have
$|\rho_k-k|<\xi_r$ for all $k\in{\mathbb N}.$ Hence, the maximum modulus principle gives
\begin{equation}\label{Delta_k-1}
|\Delta_k(k^2)| <\max_{|\rho-\rho_k|=\xi_r}\Big|\frac{\Delta(\rho^2)}{\sigma_{r,k}(\rho^2)}\Big|
=\max_{|\rho-\rho_k|=\xi_r}\Big|\Delta(\rho^2)\prod_{j\in\Omega_r(k)}\frac{j^2}{\lambda_j-\rho^2}\Big|.
\end{equation}
By virtue of representation (\ref{char1-1}) and Proposition 1, we have the estimate
$$
|\Delta(\rho^2)|\le\frac{A_{r,C}}{|\rho|+r+\xi_r}, \quad |{\rm Im}\rho|\le C,
$$
where $A_{r,C}$ depends only on $r$ and $C.$ The latter estimate holds also if $|\rho-\rho_k|=\xi_r$ for any $k\in{\mathbb N},$ because in
this case we have $|{\rm Im}\rho|\le|{\rm Im}\rho_k| +|{\rm Im}(\rho-\rho_k)|\le C=r+\xi_r.$ Furthermore, for $|\rho-\rho_k|=\xi_r$ we have
the estimate
$$
\frac1{|\rho|+r+\xi_r} \le\frac1{|\rho_k|-|\rho-\rho_k|+r+\xi_r} =\frac1{|\rho_k|+r}\le\frac1k, \quad k\in{\mathbb N}.
$$
Thus, we get the estimate
\begin{equation}\label{Delta_k-2}
|\Delta(\rho^2)|\le \frac{C_r}k, \quad |\rho-\rho_k|=\xi_r, \quad k\in{\mathbb N}.
\end{equation}
Further, we have
\begin{equation}\label{Delta_k-3}
\prod_{j\in\Omega_r(k)}j^2 <(k+5m)^{2\alpha_{r,k}}\le C_r k^{2\alpha_{r,k}}, \quad k\in{\mathbb N}.
\end{equation}
Moreover, if $|\rho-\rho_k|=\xi_r$ and $j\in\Omega_r(k),$ then
$$
|\rho-\rho_j|\ge\xi_r -|\rho_j-\rho_k| \ge\xi_r -|j-k|-\frac{|\varkappa_j|}j-\frac{|\varkappa_k|}k >\xi_r-5r=1, \quad k\ge1,
$$
$$
|\rho+\rho_j|\ge|\rho_j+\rho_k|-\xi_r \ge j+k-2r -\xi_r >2k-5r -\xi_r\qquad\qquad\qquad\quad
$$
$$
\qquad\qquad\qquad=2k-2\xi_r+1=\Big(2-\frac{2\xi_r-1}k\Big)k\ge \Big(2-\frac{2\xi_r-1}{\xi_r}\Big)k =\frac{k}{\xi_r}, \quad k\ge\xi_r.
$$
Hence, we have $|\lambda_j-\rho^2|\ge k/\xi_r$ as soon as $|\rho-\rho_k|=\xi_r,\;j\in\Omega_r(k)$ and $k\ge\xi_r,$ which along with
(\ref{Delta_k-1})--(\ref{Delta_k-3}) give (\ref{Delta_k-0}) for $k\ge\xi_r.$ Further, for $k=\overline{1,\xi_r-1}$ and $j\in\Omega_r(k)$ we
have
$$
j<k+3r\le8r, \quad |\lambda_j-\rho^2| =|\rho-\rho_j||\rho+\rho_j| \ge(10r-j-r)^2=r^2\quad {\rm for} \quad |\rho|=10r.
$$
Hence, for $k<\xi_r,$ estimate (\ref{Delta_k-0}) follows from the following rough estimate:
$$
|\Delta_k(k^2)| <\max_{|\rho|=10r}\Big|\Delta(\rho^2)\prod_{j\in\Omega_r(k)}\frac{j^2}{\lambda_j-\rho^2}\Big|
\le8^{12r-2}\max_{|\rho|=10r}|\Delta(\rho^2)|\le C_r, \quad k=\overline{1,\xi_r-1}.
$$
Thus, we arrive at (\ref{Delta_k-0}) for all $k\in{\mathbb N}.$  $\hfill\Box$

\medskip
Now we are in position to give the proof of Theorem~6.

\medskip
{\it Proof of Theorem~6}. Subtracting (\ref{tilde_Delta}) from (\ref{char1-1}), and using Parseval's equality and (\ref{Delta_k}), we get
\begin{equation}\label{1+} 
\|\hat v\|= \sqrt{\frac2\pi\sum_{k=1}^\infty k^4|\hat\Delta(k^2)|^2}, \quad \hat\Delta(k^2) =\Delta_k(k^2)\Big(\sigma_{r,k}(k^2)
-\tilde\sigma_{r,k}(k^2) +\Big(1-\frac{\tilde\Delta_k(k^2)}{\Delta_k(k^2)}\Big)\tilde\sigma_{r,k}(k^2)\Big),
\end{equation}
where
$$
\Big|1-\frac{\tilde\Delta_k(k^2)}{\Delta_k(k^2)}\Big| =\Big|1-\prod_{j\in\Theta_r(k)}\frac{\tilde\lambda_j-k^2}{\lambda_j-k^2}\Big|
=\Big|1-\exp\Big(\sum_{j\in\Theta_r(k)}\ln\Big(1-\frac{\hat\varkappa_j}{\lambda_j-k^2}\Big)\Big)\Big|,
$$
$$
\frac{|\hat\varkappa_j|}{|\lambda_j-k^2|} =\frac{2r}{|j-k|(j+k)-r} \le \frac{2r}{3r(1+1)-r}=\frac25<\frac12, \quad j\in\Theta_r(k).
$$
Thus, we get
$$
\Big|1-\frac{\tilde\Delta_k(k^2)}{\Delta_k(k^2)}\Big| \le\sum_{\nu=1}^\infty
\frac{2^\nu}{\nu!}\Big(\sum_{j\in\Theta_r(k)}\frac{|\hat\varkappa_j|}{|\lambda_j-k^2|}\Big)^\nu
=\sum_{\nu=1}^\infty\frac{(2\theta_k)^\nu}{\nu!} \le2\theta_k\exp(2\theta_k), \quad
\theta_k=\sum_{j\in\Theta_r(k)}\frac{|\hat\varkappa_j|}{|\lambda_j-k^2|}.
$$
The Cauchy--Bunyakovsky--Schwarz inequality gives $\theta_k \le\alpha_k\|\{\hat\varkappa_n\}\|_{l_2},$ $k\in{\mathbb N},$ where
$$
\alpha_k=\sqrt{\sum_{j\in\Theta_r(k)}\frac1{|\lambda_j-k^2|^2}} \le\sqrt{\sum_{j\in\Theta_r(k)}\frac1{(|j-k|(j+k)-r)^2}}
\le\sqrt{\sum_{j=1}^\infty\frac1{(3r(j+1)-r)^2}} <\frac\pi{3r\sqrt6},
$$
which along with the second equality in (\ref{1+}) and Propositions~2 and~3 give
\begin{equation}\label{hat_Delta-1}
k^2|\hat\Delta(k^2)| \le C_r\sum_{j\in\Omega_r(k)}|\hat\varkappa_j| +C_r\|\{\hat\varkappa_n\}\|_{l_2}|\varkappa_k|.
\end{equation}
Since
$$
\Big(\sum_{k=1}^n a_k\Big)^2\le n\sum_{k=1}^n a_k^2, \qquad \sum_{k=1}^\infty\sum_{j\in\Omega_r(k)} a_{k,j} =\sum_{k,j\in{\mathbb N},\,
|k-j|<6r} a_{k,j} =\sum_{j=1}^\infty\sum_{k\in\Omega_r(j)} a_{k,j}
$$
for any non-negative summands, we have
$$
\sqrt{\sum_{k=1}^\infty\Big(\sum_{j\in\Omega_r(k)}|\hat\varkappa_j|\Big)^2} \le
\sqrt{\sum_{k=1}^\infty\alpha_{r,k}\sum_{j\in\Omega_r(k)}|\hat\varkappa_j|^2}
=\sqrt{\sum_{j=1}^\infty|\hat\varkappa_j|^2\sum_{k\in\Omega_r(j)}\alpha_{r,k}} \le(6r-1)\|\{\hat\varkappa_n\}\|_{l_2},
$$
which along with the first equality in (\ref{1+}) and (\ref{hat_Delta-1}) finish the proof. $\hfill\Box$
\\

{\bf \large 6. Solution of the inverse problem. Proof of Theorems~2 and 3}\\

Solution of Inverse Problem~1 can be constructed by the following algorithm.

\medskip
{\bf Algorithm 1. }{\it Let the spectrum $\{\lambda_n\}_{n\ge1}$ of some boundary value problem ${\mathcal L}(q,M)$ along with the
potential $q(x)$ be given.

(i) In accordance with (\ref{char1}), calculate the function $v(x)$ by the formula
$$
v(x)=\frac\omega2 +\frac2\pi\sum_{k=1}^\infty \Big(k^2\Delta(k^2)+(-1)^k\frac{\omega\pi}2\Big)\cos kx, \quad \omega=\frac1\pi\int_0^\pi
q(t)\,dt,
$$
where the function $\Delta(\lambda)$ is determined by (\ref{char2});

(ii) Find the function $M(x)$ as a solution of the main equation (\ref{ma_eq}).}

\medskip
{\it Proof of Theorem~2.} Let a complex-valued function $q(x)\in L_2(0,\pi)$ and a sequence of complex numbers $\{\lambda_n\}_{n\ge1}$ of
the form (\ref{asy}) be given. According to Lemma~5, the function $\Delta(\lambda),$ constructed by formula (\ref{char2}), has the form
(\ref{char1}) with some function $v(x)\in L_2(0,\pi)$ obeying (\ref{means}). By virtue of the first part of Theorem~4, the main equation
(\ref{ma_eq}) with these functions $v(x)$ and $q(x)$ has a unique solution $M(x)\in L_{2,\pi}.$ Consider the corresponding boundary value
problem ${\mathcal L}={\mathcal L}(q,M).$ Let $\tilde\Delta(\lambda)$ be its characteristic function. Then, by virtue of Lemma~3, it has
the form
$$
\tilde \Delta(\lambda)=\frac{\sin\rho\pi}{\rho}-\omega\pi\frac{\cos\rho\pi}{2\rho^2} +\int_0^\pi \tilde v(x)\frac{\cos\rho x}{\rho^2}\,dx,
$$
where $\omega$ is determined in (\ref{asy}) and
\begin{equation}\label{6.2}
-\tilde v(\pi-x)=R(\pi,x;q,M), \quad 0<x<\pi.
\end{equation}
Comparing (\ref{6.2}) with (\ref{ma_eq}), we get $\tilde v(x)=v(x)$ a.e. on $(0,\pi)$ and, hence,
$\tilde\Delta(\lambda)\equiv\Delta(\lambda).$ Thus, the spectrum of the constructed boundary value problem ${\mathcal L}$ coincides with
the sequence $\{\lambda_n\}_{n\ge1}.$

The uniqueness of $M(x)$ follows from uniqueness of solution of the main equation (\ref{ma_eq}). $\hfill\Box$

\medskip
{\it Proof of Theorem~3.} First, let us assume that $\omega=\tilde\omega=0.$ Then, according to Theorem~6, $\|\{\varkappa_n\}\|_{l_2}\le r$
and $\|\{\tilde\varkappa_n\}\|_{l_2}\le r$ imply estimate (\ref{est_w}). In particular, Theorem~6 gives $\|v\|\le C_r$ and $\|\tilde v\|\le
C_r$ as soon as $\|\{\varkappa_n\}\|_{l_2}\le r$ and $\|\{\tilde\varkappa_n\}\|_{l_2}\le r,$ respectively. Thus, assuming also $\|q\|\le r$
and $\|\tilde q\|\le r$ and taking into account that conditions (\ref{means}) and (\ref{means_ti}) are fulfilled automatically, it falls
into assumptions of Theorem~5. Hence, estimate (\ref{est1}) holds, which along with (\ref{est_w}) give~(\ref{UFStab}).

In the general case, the shift $q_1(x):=q(x)-\omega$ transforms ${\mathcal L}(q,M)$ to the problem ${\mathcal L}(q_1,M)$ with the same
$M(x)$ and with the spectrum $\{\mu_n\}_{n\ge1},$ where $\mu_n=\lambda_n-\omega=n^2+\varkappa_n,$ $n\ge1.$ Since
$\|q_1\|\le\|q\|+\sqrt\pi|\omega|\le2\|q\|\le2r,$ according to the first part of the proof, we have the estimate $\|\hat M\|_{2,\pi}\le
C_r(\|\{\hat\varkappa_n\}\|+\|\hat q_1\|),$ while $\|\hat q_1\|\le2\|\hat q\|,$ which finishes the proof. $\hfill\Box$

\medskip
{\bf Remark 1.} Finally, for illustrating the equivalence of the group of estimates (\ref{StabYur}) to the estimate $\|\hat M_0\|\le
C\|\{\hat\lambda_n\}\|_{l_2},$ we establish the following two-sided estimate:
$$
C_1\|\hat M_0\| \le \Big(\|\hat M_0\|_{L(0,\pi)} +\|\hat M_1\|_{L(0,\pi)} +\|\hat Q\|\Big) \le C_2\|\hat M_0\|,
$$
where the functions $M_0(x),$ $M_1(x)$ and $Q(x)$ are defined in (\ref{CondYur}) and we use the designation
$\hat\gamma:=\gamma-\tilde\gamma.$ Indeed, Fubibi's theorem and the Cauchy--Bunyakovsky--Schwarz inequality imply the estimates $\|\hat
M_1\|_{L(0,\pi)}\le \|\hat M_0\|_{L(0,\pi)}\le \sqrt\pi\|\hat M_0\|,$ respectively. By virtue of Lemma~2.1 in \cite{But06_MZ}, we also have
$\|\hat M_1\|\le 2\|\hat M_0\|$ and, hence, $\|\hat Q\|\le 3\|\hat M_0\|.$ Thus, one can take $C_2=2\sqrt\pi+3.$

On the other hand, since definition of the function $\hat Q(x)$ coincides with the first relation in (\ref{eq2.2.17-1}) after putting
$z(x)=\hat Q(x)$ and $h(x)=\hat M_0(x),$ using the second relation in (\ref{eq2.2.17-1}) along with the plain fact that the mean value of
$\hat Q(x)$ on $(0,\pi)$ always vanishes, we obtain
$$
\hat M_0(x)=\hat Q(x)-\frac1{\pi-x}\int_x^\pi \hat Q(t)\,dt,
$$
which along with Lemma~2.1 in \cite{But06_MZ} give $\|\hat M_0\|\le 3\|\hat Q\|.$ Thus, it sufficient to take $C_1=1/3.$

\bigskip
{\bf Acknowledgement.} This work was supported by Grant 20-31-70005 of the Russian Foundation for Basic Research.


\begin{thebibliography}{99}

\bibitem{Mar77}
Marchenko V.A. {\it Sturm--Liouville Operators and Their Applications}, Naukova Dumka, Kiev, 1977;
English  transl., Birkh\"auser, 1986.

\bibitem{Lev84}
Levitan B.M. {\it Inverse Sturm--Liouville Problems}, Nauka, Moscow, 1984 (Russian); English transl., VNU Sci.Press, Utrecht, 1987.

\bibitem{FY01}
Freiling G. and Yurko V.A. {\it Inverse Sturm--Liouville Problems and Their Applications}, NOVA Science Publishers, New York, 2001.

\bibitem{Yur07}
Yurko V.A. {\it Introduction to the Theory of Inverse Spectral Problems}, Moscow, Fizmatlit, 2007.

\bibitem{Yur02}
Yurko V.A. {\it Method of Spectral Mappings in the Inverse Problem Theory}, Inverse and Ill-posed Problems Series. VSP, Utrecht, 2002.

\bibitem{Yur00}
Yurko V.A. {\it Inverse Spectral Problems for Differential Operators and Their Applica\-tions}, Gordon and Breach Science Publishers,
Amsterdam, 2000.

\bibitem{Mal}
Malamud M.M. {\it On some inverse problems}, Boundary Value Problems of Mathematical Physics, Kiev, 1979, 116--124.

\bibitem{Yur84}
Yurko V.A. {\it Inverse problem for integro-differential operators of the first order}, Functional Analysis, Ul'janovsk, 1984, 144--151.

\bibitem{Er}
Eremin M.S. {\it An inverse problem for a second-order integro-differential equation with a singularity}, Diff. Uravn. 24 (1988) no.2,
350--351.

\bibitem{Yur91}
Yurko V.A. {\it An inverse problem for integro-differential operators}, Mat. Zametki, 50 (1991), no.5, 134--146 (Russian); English transl.
in Math. Notes 50 (1991), no. 5--6, 1188--1197.

\bibitem{But07MM} Buterin S.A. {\it On an inverse problem for integro-differential operators of the second order}, in: Matematika.
Mekhanika, vol. 9, Saratov Univ., Saratov, 2007, pp. 8--11.

\bibitem{But07}
Buterin S.A. {\it On an inverse spectral problem for a convolution integro-differential operator}, Results Math. 50 (2007) no.3-4,
173--181.

\bibitem{Kur}
Kuryshova Ju.V. {\it Inverse spectral problem for integro-differential operators}, Mat. Zametki 81 (2007) no.6, 855--866; English transl.
in Math. Notes 81 (2007) no.6, 767--777.

\bibitem{But10}
Buterin S.A. {\it On the reconstruction of a convolution perturbation of the Sturm--Liouville operator from the spectrum}, Diff. Uravn. 46
(2010) no.1, 146--149 (Russian); English transl. in Diff. Eqns. 46 (2010) no.1, 150--154.

\bibitem{KurSh10}
Kuryshova Yu.V. and Shieh C.-T. An inverse nodal problem for integro-differential operators, J. Inverse and Ill-Posed Problems 18 (2010)
no.4, 357--369.

\bibitem{W}
Wang Y. and Wei G. {\it The uniqueness for Sturm--Liouville problems with aftereffect}, Acta Math Sci. 32A (2012) no.6, 1171--1178.

\bibitem{Yur14}
Yurko V.A. {\it An inverse spectral problems for integro-differential operators}, Far East J. Math. Sci. 92 (2014) no.2, 247--261.

\bibitem{BCh15}
Buterin S.A. and Choque Rivero A.E. {\it On inverse problem for a convolution integro-differential operator with Robin boundary
conditions}, Appl. Math. Lett. 48 (2015) 150--155.

\bibitem{BS}
Buterin S.A. and Sat M. {\it On the half inverse spectral problem for an integro-differential operator}, Inverse Problems in Science and
Engineering 25 (2017) no.10, 1508--1518.

\bibitem{BB17}
Bondarenko N. and Buterin S. {\it On recovering the Dirac operator with an integral delay from the spectrum}, Results Math. 71 (2017)
no.3-4, 1521--1529.

\bibitem{Yur17-1}
Yurko V.A. {\it Inverse spectral problems for first order integro-differential operators}, Boundary Value Problems (2017) 2017:98, 7pp.

\bibitem{BondBut-18}
Bondarenko N. and Buterin S. {\it An inverse spectral problem for integro-differential Dirac operators with general convolution kernels},
Applicable Analysis (2018), 17pp. https://doi.org/10.1080/00036811.2018.1508653

\bibitem{But18}
Buterin S.A. {\it On inverse spectral problems for first-order integro-differential operators with discontinuities}, Appl. Math. Lett. 78
(2018), 65--71.

\bibitem{But18-2} Buterin S.A. {\it Inverse spectral problem for Sturm--Liouville integro-differential operators with discontinuity
conditions}, Sovr. Mat. Fundam. Napravl. 64 (2018) no.3, 427--458; Engl. transl. in J. Math. Sci. (to appear)

\bibitem{Bon18-1}
Bondarenko N.P. {\it An inverse problem for an integro-differential operator on a star-shaped graph}, Math. Meth. Appl. Sci. 41 (2018)
no.4, 1697--1702.

\bibitem{ButVas18}
Buterin S.A. and Vasiliev S.V. {\it On uniqueness of recovering the convolution integro-differential operator from the spectrum of its
non-smooth one-dimensional perturbation}, Boundary Value Problems (2018) 2018:55, 12pp.

\bibitem{Zol18}
Zolotarev, V.A. {\it Inverse spectral problem for the operators with non-local potential}, Mathematische Nachrichten (2018), 1--21, DOI:
https://doi.org/10.1002/mana.201700029.

\bibitem{Ign18}
Ignatyev M. {\it On an inverse spectral problem for the convolution integro-differential operator of fractional order}, Results Math.
(2018) 73:34, 8pp.

\bibitem{Ign18-2}
Ignatiev M. {\it On an inverse spectral problem for one integro-differential operator of fractional order}, J. Inverse and Ill-posed Probl.
27 (2019) no.1, 17--23.

\bibitem{Bon18-2}
Bondarenko N.P. {\it An inverse problem for the integro-differential Dirac system with partial information given on the convolution
kernel}, J. Inverse Ill-Posed Probl. 27 (2019) no.2, 151--157.

\bibitem{Bon19-1}
Bondarenko N.P. {\it An inverse problem for an integro-differential pencil with polynomial eigenparameter-dependence in the boundary
condition}, Anal. Math. Phys. 9 (2019) 2227--2236.

\bibitem{Bon19-2}
Bondarenko N.P. {\it An inverse problem for an integro-differential equation with a convolution kernel dependent on the spectral
parameter}, Results Math. (2019) 74:148, 7pp.

\bibitem{Bon19-3}
Bondarenko N.P. {\it An inverse problem for the second-order integro-differential pencil}, Tamkang J. Math. 50 (2019) no.3, 223--231.

\bibitem{But19}
Buterin S.A. {\it An inverse spectral problem for Sturm--Liouville-type integro-differential operators with Robin boundary conditions},
Tamkang J. Math. 50 (2019) no.3, 207--221.

\bibitem{But20} Buterin S. {\it Uniform stability of the inverse spectral problem for a convolution integro-differential operator},
arXiv:2001.09915 [math.SP] (2020) 15pp.

\bibitem{Bor}
Borg G. {\it Eine Umkehrung der Sturm--Liouvilleschen Eigenwertaufgabe}, Acta Math. 78 (1946) 1--96.

\bibitem{BK19}
Buterin S.A. and Kuznetsova M.A. {\it On Borg's method for non-selfadjoint Sturm--Liouville operators}, Anal.
Math. Phys. 9 (2019) 2133--2150.

\bibitem{ButMal18}
Buterin S. and Malyugina M. {\it On global solvability and uniform stability of one nonlinear integral equation}, Results Math. (2018),
73:117, 19pp.

\bibitem{SavShk} Savchuk A.M. and Shkalikov A.A. {\it Inverse problems for Sturm--Liouville operators with potentials in Sobolev spaces:
Uniform stability}, Funk. Anal. i ego Pril. 44 (2010) no.4, 34--53; English transl. in Funk. Anal. Appl. 44 (2010) no.4, 270--285.

\bibitem{But06_MZ}
Buterin S.A. {\it Inverse spectral reconstruction problem for the convolution operator pertur\-bed by a one-dimensional operator},  Matem.
Zametki 80 (2006) no.5, 668--682 (Russian); English transl. in Math. Notes 80 (2006) no.5, 631--644.

\end{thebibliography}
\end{document}